\documentclass[11pt,english]{article}
\usepackage[T1]{fontenc}

\usepackage{graphicx, subfigure}
\expandafter\def\csname ver@subfig.sty\endcsname{}
\usepackage[latin9]{inputenc}
\usepackage{geometry}
\usepackage{amsmath}
\usepackage{appendix}
\usepackage{amssymb}
\usepackage{esint}
\usepackage{makecell,multirow,diagbox}
\usepackage{mathtools}
\usepackage{epstopdf}
\usepackage[latin9]{inputenc}
\usepackage{pict2e}
\usepackage{xcolor}
\usepackage{colortbl,booktabs}
\usepackage{stmaryrd}
\usepackage{esint}
\usepackage[all,cmtip]{xy}
\usepackage{mathtools}
\usepackage{float}
\usepackage{setspace}
\usepackage{authblk}
\usepackage{amsthm}
\usepackage{mathrsfs}
\usepackage{stmaryrd}
\usepackage{bbold}
\usepackage{color}
\usepackage{tikz}
\usepackage[all,cmtip]{xy}
\usepackage{algorithm}  
\usepackage{algorithmic}

\theoremstyle{plain}
\theoremstyle{plain}
\newtheorem{remark}{Remark}
\usepackage{color}
\definecolor{marin}{rgb} {0., 0.3, 0.7}
\definecolor{rouge}{rgb} {0.8, 0., 0.}
\definecolor{sepia}{rgb} {0.8, 0.5, 0.}
\usepackage[colorlinks,citecolor=sepia,linkcolor=marin,bookmarksopen,bookmarksnumbered]{hyperref}

\def\Bb{\mathbf{B}}

\def\xb{\mathbf{x}}
\def\vb{\mathbf{v}}
\def\ab{\mathbf{a}}

\def\cH{\mathcal{H}}

\def\be{\begin{equation}}
\def\ee{\end{equation}}
\def\tn{\textnormal}

\theoremstyle{definition}

\newcommand{\green}[1]{{\color[rgb]{0.1,0.6, 0.3}{#1}}}

\DeclareSymbolFont{largesymbol}{OMX}{yhex}{m}{n}
\DeclareMathAccent{\Widehat}{\mathord}{largesymbol}{"62}

\makeatletter

\usepackage{babel}
\usepackage{subfig}
\graphicspath{{Dirac/}}
           
\makeatother
\DeclareMathSizes{14}{14}{9.8}{7}

\usepackage{lipsum}

\begin{document}

\title{Canonical variables based numerical schemes for hybrid plasma models with kinetic ions and massless electrons}
\date{}
\author[1]{Yingzhe Li}
\author[1]{Florian Holderied}
\author[1]{Stefan Possanner}
\author[1,2]{Eric Sonnendr\"ucker}
\affil[1]{Max Planck Institute for Plasma Physics, Boltzmannstrasse 2, 85748 Garching, Germany}
\affil[2]{Technical University of Munich, Department of Mathematics, Boltzmannstrasse 3, 85748 Garching, Germany}

\maketitle
\begin{abstract}
We study the canonical variables based numerical schemes of a hybrid model with kinetic ions and mass-less electrons. Two equivalent formulations of the hybrid model are presented with the vector potentials in different gauges and the distribution functions depending on canonical momentum (not velocity), which constitutes a pair of canonical variables with the position variable. Particle-in-cell methods are used for the distribution functions, and the vector potentials are discretized by the finite element methods in the framework of finite element exterior calculus. Splitting methods are used for the time discretizations. It is illustrated that the second formulation is numerically superior and the schemes constructed based on the anti-symmetric bracket proposed have better conservation properties and lower noise, although the filters can be used to improve the schemes of the first formulation. 
\end{abstract}

%

\section{Introduction}

There are a lot of models proposed to describe complex physical processes in plasmas, which usually include different kinds of species and are inherently multi-scale. Among them,
hybrid models combine the advantages of kinetic and fluid models, in which some components of plasmas, such as high energy particles, are treated kinetically, while the remainder is described using fluid type equations. 
Compared to kinetic equations, hybrid models are more computationally efficient because of the fluid equations adopted and small scales ignored. Also they are more accurate than pure fluid type models, as kinetic effects of some components are included. 
There are many kinds of hybrid kinetic-fluid models for plasmas in literature in different contexts~\cite{23, 26, 31}. 

Physical laws can be equivalently represented with different sets of variables, among which canonical variables bring some advantages, such as symplectic methods~\cite{Feng, HLW} can be applied.  Transforming even part of the variables to canonical variables has some benefits, as an example of which, in this work we use a pair of canonical variables, i.e., position and canonical momentum, to describe the motion of particles. Some canonical variables based numerical schemes have been used, such as in~\cite{Leimkuhler, xuliwei}.

In magnetized kinetic plasma simulations, the distribution function  of each particle species  depends on velocity ${\mathbf v}$ or canonical momentum ${\mathbf p}$, which have the relation $m_i{\mathbf v} = {\mathbf p} - q_i{\mathbf A}$, where ${\mathbf A}$ denotes the vector potential and $m_i, q_i$ are the mass and charge of each particle species. Also the magnetic field can be represented by the magnetic field ${\mathbf B}$ or the vector potential ${\mathbf A}$ via the formula ${\mathbf B} = \nabla \times {\mathbf A}$. So for the magnetized kinetic/hybrid models, there are three kinds of formulations with different choices of unknowns, 1) $xvB$ formulation with the unknowns $f(t, {\mathbf x}, {\mathbf v})$ and ${\mathbf B}(t, {\mathbf x})$; 2) $xvA$ formulation with the unknowns $f(t, {\mathbf x}, {\mathbf v})$ and ${\mathbf A}(t, {\mathbf x})$; 3) $xpA$ formulation with the unknowns $f(t, {\mathbf x}, {\mathbf p})$ and ${\mathbf A}(t, {\mathbf x})$, where $f$ denotes the distribution function. The Poisson brackets of $xvB$ and  $xvA$ formulations, if there are,  can usually be derived from the counterparts of $xpA$ formulations by using the chain rules of the functional derivatives, and the former usually have more terms and are more complicated than the latter, which reflects the differences of coupling complexity. 

There have been some numerical simulation works based on the $xpA$ formulations in different contexts of plasma physics. 
 For the fully kinetic equations, in~\cite{Qincanonical} canonical symplectic methods are proposed based on the discretization of the canonical Poisson bracket of the Vlasov--Maxwell equations.
For the gyro-kinetic simulations~\cite{pham} with strong magnetic fields, parallel canonical momentum is used as an independent variables of the distribution functions to avoid the term of the time derivative of $A_\parallel$ in the Vlasov equation. Some recent developments include~\cite{Mishchenko, LinZ}. 
In the context of the hybrid plasma simulations, the $xvB$ and $xvA$ formulations are adopted numerically in most research works~\cite{Winske, 31, chacon1, chacon2, 26, chenliu, Porcelli, guoyongfu}, and the $xpA$ formulation has been used in~\cite{Winske} and the references therein, in which it is called the 'Hamiltonian method'.

In this work, we consider the numerical discretizations for the $xpA$ formulations of a hybrid model, in which all ions are treated kinetically, electrons are mass-less, and the electron pressure is determined by the electron density via the equation of state. 
This model has many applications in laser plasmas and astrophysics~\cite{Rambo, valentine, Schekochihin, Arzamasskiy},
and is obtained by taking quasi-neutral limit and mass-less electron limit for the more fundamental models.
For simplicity we consider the one ion species case and assume each ion has a unit charge. The commonly used $xvB$ formulation of this model after the dimensionless procedures in~\cite{LHPS} is
\begin{equation}\label{eq:hybridequations}
\begin{aligned}
&\frac{\partial f}{\partial t} = - {\mathbf v} \cdot \nabla f -({\mathbf E} + {\mathbf v} \times {\mathbf B}) \cdot \nabla_v f,\\
& \frac{\partial {\mathbf B}}{\partial t} = - \nabla \times {\mathbf E},
\end{aligned}
\end{equation}
where 
\begin{equation}\label{eq:ohmslaw}
 {\mathbf E} = -\frac{\nabla p_e}{n} -  {\mathbf u}\times {\mathbf B} +  \frac{\mathbf J}{n} \times {\mathbf B} \ (\text{Ohm's law})
\end{equation}
with $$ {\mathbf J} = \nabla \times {\mathbf B}, \quad n = \int f \mathrm{d}{\mathbf v}, \quad n{\mathbf u} = \int f{\mathbf v} \mathrm{d}{\mathbf v}.$$
Here $f$ is the distribution function of ions depending on time $t$, space ${\mathbf x}$, and velocity ${\mathbf v}$. ${\mathbf E}$ and ${\mathbf B}$ are electric and magnetic fields, respectively. The $n$ is the density of electrons, which equals to the density of ions $\int f \mathrm{d}{\mathbf v}$ because of the quasi-neutrality condition.  In the above Ohm's law, $p_e$ is the electron pressure, which is commonly taken to be $nT$ in the isothermal electron case with constant $T$ (the normalized temperature of electrons) or $Cn^\gamma$ with $\gamma \neq 1$ in the adiabatic electron case with constant $C$. In more general cases, a time evolution equation of $p_e$ should be included in the hybrid model~\eqref{eq:hybridequations} as~\cite{chacon1}.  Existing numerical methods for this formulation include current advance method~\cite{Matthews}, based on which, there is a particle-in-cell code CAMELIA~\cite{CAMELIA} and an Eulerian code~\cite{valentini}; Pegasus~\cite{Pegasus}, in which a constrained transport method is used to guarantee the divergence free property of magnetic field; mass and energy conserving particle-in-cell schemes on curvilinear meshes~\cite{chacon2}, which also conserve momentum provided the magnetic field is divergence free. Also structure-preserving particle-in-cell schemes are constructed in~\cite{LHPS} based on an anti-symmetric bracket and splitting methods, which conserve many properties at the same time, such as energy, quasi-neutrality condition, and divergence free property of magnetic field.
For more reviews about hybrid simulations, we refer the readers to the references~\cite{newreview, 31}. 
The difficulty in hybrid-kinetic codes of the $xvB$ formulation is obtaining a good approximation of  the electric field depending on particles, based on which accurate and efficient particle pushers can be constructed. 
The $xvA$ formulation can be obtained by just replacing ${\mathbf B}$ by $\nabla \times {\mathbf A}$ (${\mathbf A}$ is the vector potential) in the $xvB$ formulation~\eqref{eq:hybridequations}. 
Recently, there are some numerical schemes conserving energy, momentum, and mass proposed in~\cite{chacon1} based on the $xvA$ formulation using particle-in-cell methods combined with finite difference methods, in which an equation of electron pressure is used. As the vector potential is used, the divergence free property of the magnetic field obtained is guaranteed naturally.

In this work, the two $xpA$ formulations we investigate numerically are with different gauges of the vector potential. In the first $xpA$ formulation (or formulation I), the Weyl gauge is chosen for the vector potential, while in the second $xpA$ formulation (or formulation II), a different gauge is used due to that the electron pressure term is a gradient term in the isothermal and adiabatic electrons cases.
The characteristics of Vlasov equation in $xpA$ formulations constitute canonical Hamiltonian systems, for which symplectic methods~\cite{Feng, HLW} and discrete gradient methods~\cite{DIS} can be used to solve. However, in $xvB$ and $xvA$ formulations pushing particle is more complicated due to the complex electric field ${\mathbf E}$ depending on all the particles in the Ohm's law, which is quite common in various hybrid models~\cite{1, kaltsas}.
In~\cite{Winske} and some references therein, updating particles using the canonical momentum is introduced,
but the relation of the electron pressure term with the gauge in the isothermal and adiabatic electron cases is not mentioned. 
This idea of changing the gauge used can also be applied to the $xvA$ formulation investigated in~\cite{chacon1} with the anti-symmetric bracket presented in appendix~\ref{eq:xvAbracket} and relativistic hybrid models~\cite{rehy}.

Our discretizations follow the recent developments of structure-preserving methods for models in plasma physics with the aim of having better long term numerical behaviours. Fully discrete structure-preserving particle-in-cell methods for the Vlasov--Maxwell equations have been proposed in~\cite{DECVM, Qincanonical} using the finite difference methods in the framework of discrete exterior calculus~\cite{DEC}, and in~\cite{hevm, GEMPIC, frame} using the finite element methods in the framework of finite element exterior calculus~\cite{FEEC}. For electrostatic hybrid plasma simulations, structure preserving algorithms have been constructed in~\cite{qinxiao}.

In this work, for the two $xpA$ formulations, we discretize the vector potentials by finite element methods in the framework of finite element  exterior calculus~\cite{FEEC}, and the distribution functions are approximated as the sums of finite number of weighted particles. 
Splitting methods~\cite{HLW} are used in time to give two subsystems, for which implicit mid-point rules are used to solve.  For the formulation I with the Weyl gauge, a projector is used to deal with the electron pressure term to make it only contribute to the curl free part of the vector potential. 
Some binomial filters are applied for the  electron pressure term $T\nabla n / n$ to reduce the noise from the particle methods ($n$ is obtained by depositions of particles),  which improves significantly the schemes.
For the formulation II, discrete energy is conserved
 after the phase-space discretization and the semi-discrete system can be expressed in an anti-symmetric bracket formulation, for which the bracket splitting~\cite{2020en} is adopted. Then the midpoint rule is used for each subsystem, which is symplectic for the first subsystem and energy-conserving for the second subsystem. Our schemes have good conservation properties (such as for energy and momentum) and the magnetic field obtained is divergence free, although our schemes are not Poisson structure preserving as~\cite{DECVM, Qincanonical, hevm, GEMPIC, frame} because the bracket proposed for the continuous problem is anti-symmetric but not Poisson.

There are some connections between our recent work~\cite{LHPS} and the formulation II. The bracket proposed in~\cite{LHPS} could be derived from the bracket of formulation II using the chain rules of the functional derivatives.  
The methods constructed in this work for the formulation II have the following advantages. 1) Particles are updated by solving simple canonical Hamiltonian systems, which is good for conservation property and good long time behavior; 2) It is easier to do the implementations than~\cite{LHPS}, as only two subsystems are obtained after bracket splitting~\cite{2020en} in time, particles and fields are decoupled. 3) The algorithms obtained are more efficient, as some heavy iterations about particles and projectors for current terms in~\cite{LHPS} are avoided; 4) The property of divergence-free  of the magnetic fields is guaranteed naturally as the vector potential is used and discretized in a $H(\text{curl})$ type finite element space.

This paper is organized as follows. Two canonical momentum based formulations of the hybrid model are presented in section~\ref{sec:canonicalform}, for which phase-space and time discretizations are done in section~\ref{sec:discretization}. In section~\ref{sec:numerical},  four numerical experiments: finite grid instability, R-wave, Bernstein waves, and ion cyclotron instability are conducted to validate the codes, and comparisons are made.  In section~\ref{sec:conclusion}, we conclude the paper with a summary and an outlook to future works.

\section{Canonical momentum based formulations}\label{sec:canonicalform}
In this section, we present two equivalent formulations of the hybrid model, in which the distribution functions depend on the canonical momentum, and the vector potentials are used. For simplicity we focus on the isothermal electron case, i.e., $p_e = nT$, similar results can be obtained for the adiabatic electron case. 
As we shall see, the main differences between these two formulations are where the term $T\frac{\nabla n}{n}$ appears and how it influences the update of the particles.
The first $xpA$ formulation is obtained as follows by choosing the vector potential in Weyl gauge, this gauge is used in many plasma simulations, such as~\cite{Qincanonical, chacon1}.\\
\noindent{\bf Formulation I}: by change of unknowns as ${\mathbf B} = \nabla \times {\mathbf A}$ (${\mathbf A}$ is the Weyl gauge) and $f(t, {\mathbf x},  {\mathbf v}) = f(t, {\mathbf x},  {\mathbf p} - {\mathbf A}) =: f_{m}(t, {\mathbf x},  {\mathbf p})$, we have the equations about $f_m$ (still denoted by $f$ for convenience) and ${\mathbf A}$ from equations~\eqref{eq:hybridequations}, 
\begin{equation}\label{eq:1ee}
\begin{aligned}
& \frac{\partial f}{\partial t} = -({\mathbf p} - {\mathbf A}) \cdot \frac{\partial f}{\partial {\mathbf x}} + \left[\left( \frac{\partial {\mathbf A}}{\partial \mathbf x} \right)^\top  ({\mathbf A} - {\mathbf p})\right] \cdot \frac{\partial f}{\partial {\mathbf p}},\\
& \frac{\partial {\mathbf A}}{\partial t} = T \frac{\nabla n}{n} - \frac{\nabla \times \nabla \times {\mathbf A}}{n} \times {\nabla \times {\mathbf A}} - \frac{\int ({\mathbf A} - {\mathbf p}) f \mathrm{d}{\mathbf p}}{n} \times {\nabla \times {\mathbf A}}, \quad n = \int f \mathrm{d}{\mathbf p}.
\end{aligned}
\end{equation}
We can see that the term $T\frac{\nabla n}{n}$ influences the update of the particles through the vector potential.
The total energy of the $xpA$ formulation~\eqref{eq:1ee} is 
\begin{equation}\label{eq:1ha}
\mathcal{H} = \frac{1}{2}   \int f |{\mathbf p} - {\mathbf A}|^2   \mathrm{d}{\mathbf x} \mathrm{d}{\mathbf p} + T \int n \ln n \mathrm{d}{\mathbf x} + \frac{1}{2} \int |\nabla \times {\mathbf A}|^2   \mathrm{d}{\mathbf x}.
\end{equation}
Here the density of electrons $n$ is not regarded as an independent unknown, the reason is that when $n$ is regarded as an independent unknown, conserving the quasi-neutrality relation $n = \int f \mathrm{d}{\mathbf v}$ and the positivity of $n$ is not easy to achieve numerically, although the system~\eqref{eq:1ee} with the time evolution equation of $n$, i.e., $\partial_t n + \nabla \cdot (n{\mathbf u}) = 0$, is a Hamiltonian system with a Poisson bracket proposed in~\cite{1}. 
As $T\frac{\nabla n}{n}$ is the gradient of $T\ln n$, it has no contribution for the magnetic field $\nabla \times {\mathbf A}$. Also numerically, $n$ is obtained from particles by depositions, then $n$ and $T\frac{\nabla n}{n}$ have a lot of noise and would make the update of ${\mathbf A}$ not accurate. These inspire us to change the gauge used, and choose the gauge of ${\mathbf A}$ satisfying 
$$
{\mathbf E} = - T \frac{\nabla n}{n} - \frac{\partial {\mathbf A}}{\partial t}.
$$
Then we get the following equivalent formulation.\\
{\bf Formulation II:}
\begin{equation}\label{eq:2ee2}
\begin{aligned}
& \frac{\partial f}{\partial t} = -({\mathbf p} - {\mathbf A}) \cdot \frac{\partial f}{\partial {\mathbf x}} + \left[ T \frac{\nabla n}{n} + \left( \frac{\partial {\mathbf A}}{\partial \mathbf x} \right)^\top  ({\mathbf A} - {\mathbf p})\right] \cdot \frac{\partial f}{\partial {\mathbf p}},\\
& \frac{\partial {\mathbf A}}{\partial t} = - \frac{\nabla \times \nabla \times {\mathbf A}}{n} \times {\nabla \times {\mathbf A}} - \frac{\int ({\mathbf A} - {\mathbf p}) f \mathrm{d}{\mathbf p}}{n} \times {\nabla \times {\mathbf A}}, \quad n = \int f \mathrm{d}{\mathbf p},
\end{aligned}
\end{equation}
in which the term $T\frac{\nabla n}{n}$ appears in the Vlasov equation and has a direct interaction with particles. 
The total energy of the second $xpA$ formulation~\eqref{eq:2ee2} is 
\begin{equation}\label{eq:2ha}
\mathcal{H} = \frac{1}{2}   \int f |{\mathbf p} - {\mathbf A}|^2   \mathrm{d}{\mathbf x} \mathrm{d}{\mathbf p} + T \int n \ln n \mathrm{d}{\mathbf x} + \frac{1}{2} \int |\nabla \times {\mathbf A}|^2   \mathrm{d}{\mathbf x}.
\end{equation}
The formulation~\eqref{eq:2ee2} can be derived from the following anti-symmetric bracket~\eqref{eq:original} with energy~\eqref{eq:2ha},
\begin{equation}\label{eq:original}
\begin{aligned}
\{ \mathcal{F}, \mathcal{G} \}(f, {\mathbf A}) &=   \int f   \left[\frac{\delta \mathcal{F}}{\delta f}, \frac{\delta \mathcal{G}}{\delta f} \right]_{xp}   \mathrm{d}{\mathbf x} \mathrm{d}{\mathbf p}  - \int \frac{1}{n} \nabla \times \mathbf{A} \cdot \left(  \frac{\delta \mathcal{F}}{\delta {\mathbf A}} \times \frac{\delta \mathcal{G}}{\delta {\mathbf A}} \right) \mathrm{d}{\mathbf x},
\end{aligned}
\end{equation}
where the two terms are responsible for the derivations of the two equations in~\eqref{eq:2ee2}, respectively. 
The anti-symmetric bracket~\eqref{eq:original} is the foundation for us to use the bracket splitting~\cite{2020en}  in time for this formulation~\eqref{eq:2ee2} in the next section, which gives two sub-systems. As we shall see, the subsystem of the Vlasov equation is discretized as a canonical Hamiltonian system, for which symplectic methods~\cite{HLW} or discrete gradient methods~\cite{DIS, go} can be used for the time discretization. The term $T\frac{\nabla n}{n}$ interacts with the particles directly, and the schemes constructed are more accurate, which are presented and explained qualitatively and quantitively in the remaining parts of this work. Changing the gauge via moving out the contribution of $T\frac{\nabla n}{n}$ resembles the techniques of decomposing the parallel vector potential into the symplectic and Hamiltonian parts in~\cite{Mishchenko} and the adiabatic and non-adiabatic parts in~\cite{dongjian}.
\begin{remark}
When there is a given background magnetic field ${\mathbf B}_0$, complete equations of the above two formulations are presented in Appendix~\ref{subsec:back}. 
\end{remark}

\section{Numerical discretization}\label{sec:discretization}
In this section, we use the finite element method in the framework of finite element exterior calculus to discretize the vector potential, and particle-in-cell method to discretize the distribution function. Splitting methods are used in time for two formulations~\eqref{eq:1ee} and~\eqref{eq:2ee2}. 
Also binomial filters are introduced, which would be used in the section~\ref{sec:numerical} to reduce the noise of the particle methods of the formulation I~\eqref{eq:1ee}. Time step size is $\Delta t$, $a^n$ means the value of $a$ at $n-$th time step, and $a^{n+\frac{1}{2}}$ represents $\frac{a^n + a^{n+1}}{2}$.

\noindent{\bf Commuting diagram with B-splines}
We perform the spatial discretizations in the framework of Finite Element Exterior Calculus (FEEC).
Finite element (FE) spaces and corresponding projectors are chosen such that the following diagram commutes,
\begin{equation}
\large
\begin{aligned}
\label{diagram}
 \xymatrix{
    H^{1}({\Omega}) \ar[rr]^{{\nabla}} \ar[d]_{\Pi_0} && H(\text{curl}, {\Omega}) \ar[rr]^{{\nabla}\times} \ar[d]^{\Pi_1} && H(\text{div}, {\Omega}) \ar[rr]^{{\nabla}\cdot} \ar[d]^{\Pi_2}  && L^2({\Omega}) \ar[d]^{\Pi_3} \\
      V_0 \ar[rr]^{{\nabla}}  && V_1 \ar[rr]^{{\nabla}\times}  && V_2  \ar[rr]^{{\nabla}\cdot} && V_3 
}
\end{aligned}
\end{equation}
where $V_0, V_1, V_2$ and $V_3$ are finite element spaces in which fields (proxies of $p$-forms, $p\leq 3$) are discretized in. The projectors $\Pi_n$ are based on inter-/histopolation at/between Greville points of the B-splines which span the FE spaces. For details we refer to \cite{3} which uses exactly the same basis functions and projectors. The FE spaces are written as
\begin{alignat*}{2}
& V_0 := \text{span} \{\Lambda^0_i | 0 \le i < N_0  \}, &&\ni {a}^0_h(t, {\mathbf x}) = \sum_{i=0}^{N_0-1} a_i(t)\Lambda^0_i({\mathbf x}) \\
& V_1 := \text{span}  \left\{
\left(\begin{matrix}
  \Lambda^1_{1,i} \\
    0\\
      0
\end{matrix} \right),
\left(\begin{matrix}
  0\\
     \Lambda^1_{2,i} \\
      0
\end{matrix} \right),
\left(\begin{matrix}
  0 \\
    0\\
       \Lambda^1_{3,i} 
\end{matrix} \right)
\Bigg| 
\begin{matrix}
  0 \le i < N^{1}_1\\
    0 \le i < N^{1}_2 \\
      0 \le i < N^{1}_3
\end{matrix} 
\right\}
&&\ni {\mathbf a}^1_h(t, {\mathbf x}) = \sum_{\mu = 1}^3 \sum_{i=1}^{N_\mu^1-1} a_{\mu,i}(t) \Lambda^1_{\mu, i}({\mathbf x}) {\mathbf e}_\mu,\\
& V_2 := \text{span}  \left\{
\left(\begin{matrix}
  \Lambda^2_{1,i} \\
    0\\
      0
\end{matrix} \right),
\left(\begin{matrix}
  0\\
     \Lambda^2_{2,i} \\
      0
\end{matrix} \right),
\left(\begin{matrix}
  0 \\
    0\\
       \Lambda^2_{3,i} 
\end{matrix} \right)
\Bigg| 
\begin{matrix}
  0 \le i < N^{2}_1\\
    0 \le i < N^{2}_2 \\
      0 \le i < N^{2}_3
\end{matrix} 
\right\}
&&\ni {\mathbf a}^2_h(t, {\mathbf x}) = \sum_{\mu= 1}^3 \sum_{i=1}^{N_\mu^1-1} a_{\mu,i}(t) \Lambda^2_{\mu, i}({\mathbf x}) {\mathbf e}_\mu,\\
& V_3 := \text{span} \{\Lambda^3_i | 0 \le i < N_3  \},  &&\ni {a}^3_h(t, {\mathbf x}) = \sum_{i=0}^{N_3-1} a_i(t)\Lambda^3_i({\mathbf x}) \,.
\end{alignat*}
Here, the functions $\Lambda_i^n : \Omega \to \mathbb R$ are tensor products of uni-variate B-splines of different degree, as described in \cite{3, GEMPIC}, and ${\mathbf e}_1 = (1,0,0)^\top$, ${\mathbf e}_2 = (0,1,0)^\top$, ${\mathbf e}_3 = (0,0,1)^\top$. The dimensions are
\be
 \tn{dim} V_0 = N_0\,,\quad \tn{dim} V_1 = N_1 =  \sum_{d=1}^3 N^1_d\,, \quad \tn{dim} V_2 = N_2 =  \sum_{d=1}^3 N^2_d\,,\quad \tn{dim} V_3 = N_3\,.
\ee
To simplify the notation, we stack the FE coefficients $a_i$ and basis functions in column vectors, e.g. ${\mathbf a} := (a_i)_{0 \le i  <N_0} \in \mathbb{R}^{N_0}$, and ${{\mathbf \Lambda}^0} := (\Lambda^0_i)_{0\le  i  < N_0} \in \mathbb{R}^{N_0}$. Spline functions can then be compactly written as  
\begin{equation*}
\begin{aligned}
& {a}_h^0 = {\mathbf a}^{\top} {\mathbf  \Lambda}^0,\\
& ({\mathbf a}_h^1)^\top = (\underbrace{a_{1,0}, \cdots, a_{1,N^1_1-1}}_{=: {\mathbf a}_1^\top}, \underbrace{a_{2,0}, \cdots, a_{2,N^1_2-1}}_{=: {\mathbf a}_2^\top}, \underbrace{a_{3,0}, \cdots, a_{3,N^1_3-1}}_{=: {\mathbf a}_3^\top}) 
\left(\begin{matrix}
  {\mathbf \Lambda}^1_{1} & 0 & 0 \\
    0 & {\mathbf \Lambda}^1_{2}  & 0 \\
      0 & 0 & {\mathbf \Lambda}^1_{3}
\end{matrix} \right) =: {\mathbf a}^\top \mathbb{\Lambda}^1,\\
& ({\mathbf a}_h^2)^\top = (\underbrace{a_{1,0}, \cdots, a_{1,N^2_1-1}}_{=: {\mathbf a}_1^\top}, \underbrace{a_{2,0}, \cdots, a_{2,N^2_2-1}}_{=: {\mathbf a}_2^\top}, \underbrace{a_{3,0}, \cdots, a_{3,N^2_3-1}}_{=: {\mathbf a}_3^\top}) 
\left(\begin{matrix}
  {\mathbf \Lambda}^1_{1} & 0 & 0 \\
    0 & {\mathbf \Lambda}^1_{2}  & 0 \\
      0 & 0 & {\mathbf \Lambda}^1_{3}
\end{matrix} \right) =: {\mathbf a}^\top \mathbb{\Lambda}^2,\\
& {a}_h^3 = {\mathbf a}^{\top} {\mathbf  \Lambda}^3,
\end{aligned}
\end{equation*}
where $\mathbb{\Lambda}^1 \in \mathbb{R}^{N_1 \times 3}$ and $\mathbb{\Lambda}^2 \in \mathbb{R}^{N_2 \times 3}$.
In this setting the discrete representations of the exterior derivatives can be written as matrices solely acting on finite element coefficients, 
$$
V_1\ni {\nabla} {a}^0_{h} = ({\mathbb{G} {\mathbf {\mathbf a}}})^\top \mathbb{\Lambda}^1, \quad V_2 \ni  {\nabla} \times {\mathbf a}_h^1 = ({\mathbb{C} {\mathbf a}})^\top \mathbb{\Lambda}^2\,,\quad V_3 \ni  \nabla \cdot  \ab^2_h = (\mathbb D \ab)^\top \mathbb \Lambda^3\,.
$$
where $\mathbb{G} \in \mathbb{R}^{N_1 \times N_0}$,  $\mathbb{C} \in \mathbb{R}^{N_2 \times N_1}$ and $\mathbb{D} \in \mathbb{R}^{N_3 \times N_2}$ are sparse and contain only zeros and ones.
Finally, the (symmetric) mass matrices corresponding to the discrete spaces $V_0$-$V_3$ follow from the $L^2$-inner products of basis functions,
\begin{alignat}{2}
&\mathbb{M}_0 := \int {\mathbf \Lambda^0} ({\mathbf \Lambda^0})^\top \,\tn{d}^3{\mathbf x} &&\in \mathbb R^{N_0 \times N_0}\,,
\\
&\mathbb{M}_1 := \int {\mathbb{\Lambda}^1} ({\mathbb{\Lambda}^1})^\top \,\tn{d}^3{\mathbf x} &&\in \mathbb R^{N_1 \times N_1}\,,
\\
&\mathbb{M}_2 := \int {\mathbb{\Lambda}^2}  ({\mathbb{\Lambda}^2})^\top \, \tn{d}^3{\mathbf x} &&\in \mathbb R^{N_2 \times N_2}\,,
\\
&\mathbb{M}_3 := \int {\mathbf \Lambda^3} ({\mathbf \Lambda^3})^\top \, \tn{d}^3{\mathbf x} &&\in \mathbb R^{N_3 \times N_3}\,.
\end{alignat}
These mass matrices are sparse because of the compact supports of B-splines.

\noindent{\bf Particle-in-cell methods}
The distribution function is discretized by particle-in-cell methods with delta functions, i.e., 
\begin{equation}
f(t, {\mathbf x}, {\mathbf p}) \approx f_h(t, {\mathbf x}, {\mathbf p}) = \sum_{k=1}^K w_k \delta({\mathbf x} - {\mathbf x}_k)  \delta({\mathbf p} - {\mathbf p}_k), 
\end{equation}
or smoothed delta functions, i.e., 
\begin{equation}
f(t, {\mathbf x}, {\mathbf p}) \approx f_h(t, {\mathbf x}, {\mathbf p}) = \sum_{k=1}^K w_k S({\mathbf x} - {\mathbf x}_k)  \delta({\mathbf p} - {\mathbf p}_k), 
\end{equation}
where $K$ is the total particle number, and constant $w_k, 1 \le k \le K$ represents the weight of $k$-th particle.
Smoothed delta function $S$ is defined as 
\begin{equation}
S(\mathbf  x) = \frac{1}{h_1h_2h_3} S_{k_1}\left(\frac{x_1}{h_1}\right)S_{k_2}\left(\frac{x_2}{h_2}\right)S_{k_3}\left(\frac{x_3}{h_3}\right), 
\end{equation}
where $S_k$ is defined as 
$$
S_0(x) := \mathbb{1}_{[-\frac{1}{2}, \frac{1}{2}]}, \quad S_k(x) = S_0 \star S_{k-1} = \int_{-\frac{1}{2}}^{\frac{1}{2}} S_{k-1}(x - y) \mathrm{d}{y}.
$$
Then we know that the localized support as $S({\mathbf x})$ is 
$$\text{supp}(S) = \left[-h_1\frac{k_1+1}{2}, h_1\frac{k_1+1}{2} \right] \times  \left[-h_2\frac{k_2+1}{2}, h_2\frac{k_2+1}{2} \right] \times  \left[-h_3\frac{k_3+1}{2}, h_3\frac{k_3+1}{2}\right],
$$  
where the $h_1, h_2, h_3$ are chosen as the cell sizes of the fields' discretization. 

\noindent{\bf Discrete Hamiltonian}
The vector potential ${\mathbf A}$ is regarded as a 1-form and discretized in the finite element space $V_1$, then we have
$$
{\mathbf A} \approx {\mathbf A}_h = ({\mathbb \Lambda}^1)^\top {\mathbf a}.
$$
The density of electron is approximated as $n \approx n_h = \sum_{k=1}^K w_k S({\mathbf x} - {\mathbf x}_k)$.
Then we have the discrete Hamiltonian 
\begin{equation}
\label{eq:dishhh}
\begin{aligned}
H({\mathbf X}, {\mathbf P}, {\mathbf a}) & = \frac{1}{2}\sum_{k=1}^K w_k |{\mathbf p}_k|^2 + \frac{1}{2}\sum_{k=1}^K w_k |{\mathbf A}_h({\mathbf x}_k)|^2 - \sum_{k=1}^K w_k {\mathbf p}_k \cdot {\mathbf A}_h({\mathbf x}_k) + H_{e} + \frac{1}{2}{\mathbf a}^\top \mathbb{C}^\top \mathbb{M}_2 \mathbb{C} {\mathbf a},\\
H_{e} & = T \int {n}_h \ln {n}_h \mathrm{d}{\mathbf x} \approx T \sum_j w_j \left( \sum_{k} w_k S({\mathbf x}_j - {\mathbf x}_k) \right) \ln \left( \sum_{k} w_k S({\mathbf x}_j - {\mathbf x}_k) \right),
\end{aligned}
\end{equation}
where $H_e$ is the discrete electron thermal energy, the $x_j$ and $w_j$ are quadrature points and weights. 
The energy~\eqref{eq:dishhh} can be written in a more compact way by defining suitable matrices and vectors,
\begin{equation}\label{eq:dish}
H = \frac{1}{2} {\mathbf P}^\top \mathbb{W} {\mathbf P} + \frac{1}{2}{\mathbf a}^\top {\mathbb P}_1^\top \mathbb{W}  {\mathbb P}_1 {\mathbf a}   -   {\mathbf P}^\top \mathbb{W}  \mathbb{P}_1 {\mathbf a} + H_{e} + \frac{1}{2} {\mathbf a}^\top {\mathbb C}^\top \mathbb{M}_2 {\mathbb C} {\mathbf a},
\end{equation}
where 
\begin{equation}
\label{notations}
\begin{aligned}
&{\mathbf X} := (x_{1,1},  \cdots, x_{K,1}, x_{1,2}, \cdots, x_{K,2}, x_{1,3}, \cdots, x_{K,3})^\top \quad && \in \mathbb{R}^{3K}\,,
\\[2mm]
& {\mathbf P} :=   (p_{1,1},  \cdots, p_{K,1}, p_{1,2}, \cdots, p_{K,2}, p_{1,3}, \cdots, p_{K,3})^\top \quad && \in \mathbb{R}^{3K}\,,
\\[2mm]
&\mathbb{P}^n_\mu({\mathbf X}) := (\Lambda^n_{\mu,i}({\mathbf x}_k))_{0\le i < N^n_\mu, 1 \le k \le K} \quad (n \in \{1, 2\}, \mu \in \{ 1, 2, 3\}) \quad && \in \mathbb{R}^{N_\mu^n \times K}\,,
\\[2mm]
& \mathbb{P}_n({\mathbf X}) := \text{diag}(\mathbb{P}^n_1, \mathbb{P}^n_2, \mathbb{P}^n_3), n \in \{1, 2\} \quad && \in \mathbb{R}^{N^n \times 3K}\,,
\\[2mm]
& \mathbb{W} :=  \mathbb{I}_3 \otimes \text{diag}(w_1, \cdots,  w_{K}) \quad && \in \mathbb{R}^{3K \times 3 K}.
\end{aligned}
\end{equation}

\subsection{Phase-space discretization: Formulation I}
For the formulation I~\eqref{eq:1ee}, we use splitting methods in time, and get the following two subsystems. Midpoint rule is used for time discretization, and a local projector is used for  term $T\frac{\nabla n}{n}$ to make it live in finite element space $V_1$, and some binomial filters are used to reduce the particle noise. 

\noindent{\bf The first subsystem} is
\begin{equation}\label{eq:ee}
\begin{aligned}
& \frac{\partial f}{\partial t} = -({\mathbf p} - {\mathbf A}) \cdot \frac{\partial f}{\partial {\mathbf x}} + \left[\left( \frac{\partial {\mathbf A}}{\partial \mathbf x} \right)^\top  ({\mathbf A} - {\mathbf p})\right] \cdot\frac{\partial f}{\partial {\mathbf p}},\\
& \frac{\partial {\mathbf A}}{\partial t} = T \frac{\nabla n}{n}.
\end{aligned}
\end{equation}
As we use particle-in-cell methods to discretize $f$, we have the following equations for $k$-th particle,
$$
\dot{\mathbf x}_k =  {\mathbf p}_k - {\mathbf A}_h({\mathbf x}_k), \quad \dot{\mathbf p}_k = - \left[\left( \frac{\partial {\mathbf A}_h}{\partial \mathbf x}({\mathbf x}_k) \right)^\top  ({\mathbf A}_h({\mathbf x}_k) - {\mathbf p}_k)\right], \quad 1 \le k \le K.
$$ 
As the fields ${\mathbf A}$ is regarded as a one form and discretized in finite element space $V_1$, i.e., ${\mathbf A} \approx {\mathbf A}_h = ({\mathbb \Lambda}^1)^\top {\mathbf a}$, $T \frac{\nabla n}{n}$ should also be discretized in $V_1$, where $n \approx \sum_k w_k S({\mathbf x} - {\mathbf x}_k)$, $1 \le k \le K$. To make the discretization of $T \frac{\nabla n}{n}$ live in $V_1$, a local projector is used, i.e., 
\begin{equation}\label{eq:nfilpro}
T \frac{\nabla n}{n} = T \nabla \ln n \approx T \nabla \Pi_0 \left( \ln  \left( \sum_k w_k S({\mathbf x} - {\mathbf x}_k) \right) \right) = T \Pi_1 \left(\frac{\nabla \sum_k w_k S({\mathbf x} - {\mathbf x}_k)}{ \sum_k w_k S({\mathbf x} - {\mathbf x}_k)}\right),
\end{equation}
where the last equality comes from the commuting property of the diagram~\eqref{diagram}.
The density of electrons $n_h = \sum_k w_k S({\mathbf x} - {\mathbf x}_k)$ has a lot of noise, which would make the update of ${\mathbf A}$ not stable. To solve this issue, we apply the binomial filters~\cite{Vay} to smooth the above density obtained from particles, i.e., 
$$
  T \nabla \Pi_0 \left( \ln  \left( \sum_k w_k S({\mathbf x} - {\mathbf x}_k) \right) \right) \approx  T \nabla \Pi_0 \left( \ln  \left(F \left(\sum_k w_k S({\mathbf x} - {\mathbf x}_k)\right) \right) \right),
$$
where the $F$ is the binomial filter operator.
Time discretization is done using mid-point rule,
\begin{equation}
\label{eq:xpamid}
\begin{aligned}
& \frac{{\mathbf x}_k^{n+1} - {\mathbf x}_k^{n}}{\Delta t} = {\mathbf p}_k^{n+\frac{1}{2}} - {\mathbf A}_h^{n+\frac{1}{2}}({\mathbf x}_k^{n+\frac{1}{2}}),\\
& \frac{{\mathbf p}_k^{n+1} - {\mathbf p}_k^{n}}{\Delta t} = - \left( \frac{\partial {\mathbf A}_h^{n+\frac{1}{2}}}{\partial \mathbf x}({\mathbf x}_k^{n+\frac{1}{2}}) \right)^\top  ({\mathbf A}^{n+\frac{1}{2}}_h({\mathbf x}_k^{n+\frac{1}{2}}) -  {\mathbf p}_k^{n+\frac{1}{2}}),\\
& \frac{{\mathbf a}^{n+1} - {\mathbf a}^{n}}{\Delta t} = T \mathbb{G} \tilde{\Pi}_0 \left( \ln \left(F \left(\sum_{k=1}^K w_k S({\mathbf x} - {\mathbf x}_k^{n+\frac{1}{2}} )\right)\right)\right),
\end{aligned}
\end{equation}
where $\tilde{\Pi}_0$ gives of the finite element coefficients obtained from $\Pi_0$.\\
We use the following Picard iteration to solve~\eqref{eq:xpamid}, 
\begin{equation*}
\begin{aligned}
& \frac{{\mathbf x}_k^{n+1, i+1} - {\mathbf x}_k^{n}}{\Delta t} = {\mathbf p}_k^{n+\frac{1}{2}, i} - {\mathbf A}_h^{n+\frac{1}{2},i}({\mathbf x}_k^{n+\frac{1}{2},i}),\\
& \frac{{\mathbf p}_k^{n+1,i+1} - {\mathbf p}_k^{n}}{\Delta t} = - \left( \frac{\partial {\mathbf A}_h^{n+\frac{1}{2},i}}{\partial \mathbf x}({\mathbf x}_k^{n+\frac{1}{2},i}) \right)^\top  ({\mathbf A}^{n+\frac{1}{2},i}_h({\mathbf x}_k^{n+\frac{1}{2},i}) -  {\mathbf p}_k^{n+\frac{1}{2},i}),\\
& \frac{{\mathbf a}^{n+1,i+1} - {\mathbf a}^{n}}{\Delta t} = T \mathbb{G} \tilde{\Pi}_0 \left( \ln \left(F \left(\sum_{k=1}^K w_k S({\mathbf x} - {\mathbf x}_k^{n+\frac{1}{2},i} )\right)\right)\right),
\end{aligned}
\end{equation*}
where $i$ is the iteration index, and 
$${\mathbf x}_k^{n+\frac{1}{2}, i} = ({\mathbf x}_k^n + {\mathbf x}_k^{n+1,i})/2, \quad {\mathbf p}_k^{n+\frac{1}{2}, i} = ({\mathbf p}_k^n + {\mathbf p}_k^{n+1,i})/2, \quad {\mathbf A}_h^{n+\frac{1}{2},i} = (  {\mathbf A}_h^{n} + {\mathbf A}_h^{n+1,i}    )/2.$$
We denote the solution map of this subsystem as $\Phi_{xpa}^{\Delta t}$,

\begin{remark}
In the continuous case, we have $\nabla \times (T \frac{\nabla n }{n}) = 0$, which means that $T \frac{\nabla n}{n}$ only contributes to the curl free part of ${\mathbf A}$. By~\eqref{eq:nfilpro}, the discretization $ T \nabla \Pi_0 \left( \ln  \left( F(\sum_k w_k S({\mathbf x} - {\mathbf x}_k)) \right) \right)$ is also only related with the curl free part of ${\mathbf A}$, which is consistent with the continuous case.
\end{remark}

\noindent{\bf The second subsystem} is
\begin{equation}
\label{eq:iseond}
\begin{aligned}
& \frac{\partial f}{\partial t} = 0,\\
& \frac{\partial {\mathbf A}}{\partial t} = -\frac{1}{n} {\nabla \times \nabla \times {\mathbf A}} \times {\nabla \times {\mathbf A}}- \frac{1}{n} {\int ({\mathbf A} - {\mathbf p}) f \mathrm{d}{\mathbf p}} \times {\nabla \times {\mathbf A}}.
\end{aligned}
\end{equation}
We discretize the vector potential by finite element method in weak formulation as follows, and show the energy conservation property by writing the discretization in the form of $\frac{\partial {\mathbf a}}{\partial t} = \mathbb{A}\nabla_{\mathbf a}H$, where matrix $\mathbb{A}$ is anti-symmetric. In this sub-step, $n$ is approximated as $n_h = \sum_{k=1}^K w_k S({\mathbf x} - {\mathbf x}_k)$.\\
Multiplying a test function $\mathbf{C} = (\mathbb{\Lambda}^1)^\top {\mathbf c} \in V_1$  gives,
\begin{align*}
& {\mathbf c}^\top \mathbb{M}_1 \frac{\partial {\mathbf a}}{\partial t} = \int \frac{\partial {\mathbf A}_h}{\partial t} \cdot {\mathbf C} \mathrm{d}{\mathbf x}  \\
& =  \int \left( - \frac{1}{n_h} \left( {\nabla} \times ( {\nabla} \times {\mathbf A}_h) \right) \times \left( {\nabla} \times  {\mathbf A}_h\right)  -  \frac{1}{n_h}  \int ({\mathbf A}_h - {\mathbf p}) f \mathrm{d}{\mathbf p} \times \left( {\nabla} \times  {\mathbf A}_h \right)\right)^\top  {\mathbf C} \mathrm{d} {\mathbf x}\\
& = - \underbrace{\int \left( \frac{1}{n_h}  \left( {\nabla} \times ({\nabla} \times {\mathbf A}_h) \right) \times \left( {\nabla} \times {\mathbf A}_h\right)  \right)^\top  {\mathbf C} \mathrm{d} {\mathbf x}}_{\text{term 1}} - \underbrace{ \int  \left( \frac{1}{n_h} \int ({\mathbf A}_h - {\mathbf p}) f \mathrm{d}{\mathbf p} \times \left({\nabla} \times {\mathbf A}_h \right)\right)^\top {\mathbf C} \mathrm{d} {\mathbf x}}_{\text{term 2}}.
\end{align*}
We project $\nabla \times \nabla \times {\mathbf A}_h $ into $V_1$ space by a $L^2$ projection, and have 
\begin{align*}
\text{term 1} &=  \int \left( - \frac{1}{n_h} \left(     {\nabla} \times (  {\nabla} \times {\mathbf A}_h) \right) \times \left( {\nabla} \times  {\mathbf A}_h \right)  \right) \cdot  {\mathbf C}    \mathrm{d} {\mathbf x}\\
& \approx  \int \left(  \frac{1}{n_h}    {\nabla} \times  {\mathbf A}_h   \right) \cdot  \left(  \Pi_{L^2}\left(     {\nabla} \times (  {\nabla} \times {\mathbf A}_h)  \right) \times {\mathbf C} \right)  \mathrm{d} {\mathbf x} \\
& = {\mathbf c}^\top  \mathbb{F}({\mathbf a}) \mathbb{M}_1^{-1} \left( \mathbb{C}^\top \mathbb{M}_2 \mathbb{C} {\mathbf a} \right),
\end{align*}
where 
\begin{align*}
& \mathbb{F}({\mathbf a})_{ij} = \int \left( \nabla \times {\mathbf A}_h \right) \cdot   (\Lambda^1_j  \times  \Lambda^1_i) \frac{1}{n_h} \mathrm{d} {\mathbf x}, \quad  \Pi_{L^2}\left(     {\nabla} \times (  {\nabla} \times {\mathbf A}_h) \right) = (\mathbb{\Lambda}^1)^\top \mathbb{M}_1^{-1} \left( \mathbb{C}^\top \mathbb{M}_2 \mathbb{C} {\mathbf a} \right).
\end{align*}
We  project $\int ({\mathbf A}_h - {\mathbf p}) f_h \mathrm{d}{\mathbf p}$ in term 2 into $V_1$ firstly by a $L^2$ projection, then do the similar calculations for term 1 and get 
\begin{align*}
\text{term 2} & = \int - \left( \frac{1}{n_h}  \left(\int ({\mathbf A}_h - {\mathbf p}) f_h \mathrm{d}{\mathbf p} \right) \times (\nabla \times {\mathbf A}_h) \right)^\top  {\mathbf C} \mathrm{d} {\mathbf x}\\
& \approx \int \left(    {\nabla} \times  {\mathbf A}_h   \right) \cdot  \left(  \frac{1}{n_h}   \Pi_{L^2} \left(\int ({\mathbf A}_h - {\mathbf p}) f_h \mathrm{d}{\mathbf p} \right) \times {\mathbf C} \right)   \mathrm{d} {\mathbf x} \\
& = {\mathbf c}^\top   \mathbb{F}({\mathbf a}) \mathbb{M}_1^{-1} \left(  {\mathbb P}_1^\top \mathbb{W}  {\mathbb P}_1 {\mathbf a}  -   \mathbb{P}_1^\top \mathbb{W} {\mathbf P} \right),
\end{align*}
where 
\begin{align*}
& \Pi_{L^2} \left(\int ({\mathbf A}_h - {\mathbf p}) f_h \mathrm{d}{\mathbf p} \right)
 =  \mathbb{\Lambda}^{1,\top} \mathbb{M}_1^{-1} \int  ({\mathbf A}_h - {\mathbf p}) f_h \mathrm{d}{\mathbf p} \mathbb\Lambda^1  \mathrm{d} {\mathbf x}
 =  \mathbb{\Lambda}^{1,\top} \mathbb{M}_1^{-1}  \left(  {\mathbb P}_1^\top \mathbb{W}  {\mathbb P}_1 {\mathbf a}  -   \mathbb{P}_1^\top \mathbb{W} {\mathbf P}  \right).
\end{align*}
Then Term 1 + Term 2  gives
\begin{equation}\label{eq:anti-sub2}
\frac{\partial {\mathbf a}}{\partial t} = \mathbb{M}_1^{-1}    \mathbb{F}({\mathbf a}) \mathbb{M}_1^{-1} \nabla_{\mathbf a} H,
\end{equation}
where matrix $ \mathbb{M}_1^{-1}    \mathbb{F}({\mathbf a}) \mathbb{M}_1^{-1} $ is anti-symmetric, and thus energy is conserved.
\begin{remark}
From term 1 and term 2, we get the same  matrix $\mathbb{F}({\mathbf a}) \mathbb{M}_1^{-1}$, which is quite important to get the above formulation~\eqref{eq:anti-sub2}, as $\nabla_{\mathbf a} H$ contains two terms, which are distributed in term 1 and term 2 respectively,
$$
\nabla_{\mathbf a}H = \underbrace{ \mathbb{C}^\top \mathbb{M}_2 \mathbb{C} {\mathbf a}}_{included\ in\ term\ 1} +  \underbrace{ {\mathbb P}_1^\top \mathbb{W}  {\mathbb P}_1 \frac{{\mathbf a}^n + {\mathbf a}^{n+1}}{2} -   \mathbb{P}_1^\top \mathbb{W} {\mathbf P}^n}_{included\ in\ term\ 2}.
$$
Also in the continuous PDE level,  $ \frac{\delta H}{\delta {\mathbf A}}  $ is the sum of two terms, i.e., 
$$
 \frac{\delta H}{\delta {\mathbf A}}
 = \nabla \times \nabla \times {\mathbf A} + \int \left( {\mathbf A} - {\mathbf p} \right) f \mathrm{d}{\mathbf p},
 $$
which is different from the cases in~\cite{GEMPIC, 3}. 
\end{remark}

By using mid-point rule in time, we have the following energy-conserving scheme,
\begin{equation}
\label{eq:picardaa}
\frac{{\mathbf a}^{n+1} -  {\mathbf a}^{n}  }{\Delta t} = \mathbb{M}_1^{-1} \mathbb{F}\left({\mathbf a}^{n+\frac{1}{2}}\right) \mathbb{M}_1^{-1}   \left( {\mathbb P}_1^\top \mathbb{W}  {\mathbb P}_1 \frac{{\mathbf a}^n + {\mathbf a}^{n+1}}{2} -   \mathbb{P}_1^\top \mathbb{W} {\mathbf P}^n +  {\mathbb C}^\top \mathbb{M}_2 {\mathbb C} \frac{{\mathbf a}^n + {\mathbf a}^{n+1}}{2} \right),
\end{equation}
i.e., 
\begin{equation*}
\tiny
\begin{aligned}
\left( \mathbb{M}_1 - \frac{\Delta t}{2} \left(\mathbb{F}\left({\mathbf a}^{n+\frac{1}{2}}\right)  \mathbb{M}_1^{-1} (  {\mathbb P}_1^\top \mathbb{W}  {\mathbb P}_1 +   {\mathbb C}^\top \mathbb{M}_2 {\mathbb C})  \right) \right) {\mathbf a}^{n+1} = \mathbb{M}_1{\mathbf a}^n + \frac{\Delta t}{2} \left(\mathbb{F}\left({\mathbf a}^{n+\frac{1}{2}}\right)   \mathbb{M}_1^{-1}  ({\mathbb P}_1^\top \mathbb{W}  {\mathbb P}_1 +   {\mathbb C}^\top \mathbb{M}_2 {\mathbb C})  \right) {\mathbf a}^{n} - \mathbb{F}\left({\mathbf a}^{n+\frac{1}{2}} \right) \mathbb{M}_1^{-1} {\mathbb P}_1^\top \mathbb{W}  {\mathbf P}^n.
\end{aligned}
\end{equation*}
Picard iteration and GMRES method can be used to solve this nonlinear system~\cite{iter}.
We denote the solution map of this subsystem as $\Phi_{a}^{\Delta t}$.
 
\subsection{Phase-space discretization: Formulation II}
For this formulation~\eqref{eq:2ee2}, we use bracket splitting~\cite{2020en}, instead of the Hamiltonian  splitting~\cite{ha1, ha2, ha3}, based on the bracket~\eqref{eq:original} and get two subsystems,
\begin{equation}
\begin{aligned}
& \dot{\mathcal Z} =  \{ \mathcal{Z}, \mathcal{H} \}_1,\quad \dot{\mathcal Z} =  \{ \mathcal{Z}, \mathcal{H} \}_2, \quad \mathcal{Z} = (f, {\mathbf A}),
\end{aligned}
\end{equation}
where 
$$ \{ \mathcal{F}, \mathcal{G} \}_1 =  \int f   \left[\frac{\delta \mathcal{F}}{\delta f}, \frac{\delta \mathcal{G}}{\delta f} \right]_{xp}   \mathrm{d}{\mathbf x} \mathrm{d}{\mathbf p}, \quad \{ \mathcal{F}, \mathcal{G} \}_2=   - \int \frac{1}{n} \nabla \times \mathbf{A} \cdot \left(  \frac{\delta \mathcal{F}}{\delta {\mathbf A}} \times \frac{\delta \mathcal{G}}{\delta {\mathbf A}} \right) \mathrm{d}{\mathbf x}.$$
The second subsystem $\dot{\mathcal Z} =  \{ \mathcal{Z}, \mathcal{H} \}_2$ is the same as in~\eqref{eq:iseond}, and discretizing this sub-bracket $ \{ \cdot, \cdot \}_2$ as~\cite{LHPS} by approximating the functional derivative $\frac{\delta \mathcal{F}}{\delta \mathcal{\mathbf A}}$ as $(\mathbb{\Lambda})^{1,\top} \mathbb{M}_1^{-1} \nabla_{\mathbf a} F$ gives the same discretization of~\eqref{eq:anti-sub2}.

The first sub-system $\dot{\mathcal Z} =  \{ \mathcal{Z}, \mathcal{H} \}_1$ is,
\begin{equation}\label{eq:ee2}
\begin{aligned}
& \frac{\partial f}{\partial t} = -({\mathbf p} - {\mathbf A}) \cdot \frac{\partial f}{\partial {\mathbf x}} + \left[T  \frac{\nabla n}{n} + \left( \frac{\partial {\mathbf A}}{\partial \mathbf x} \right)^\top  ({\mathbf A} - {\mathbf p})\right] \cdot \frac{\partial f}{\partial {\mathbf p}}, \quad  \frac{\partial {\mathbf A}}{\partial t} = 0,
\end{aligned}
\end{equation}
which is a Hamiltonian system.
By the particle discretization of the above $\{ \cdot, \cdot \}_1$ as~\cite{GEMPIC}, we get the discrete bracket
$$
\{F, G\}_{1,h} = \sum_{k=1}^{K} \frac{1}{w_k} \left( \nabla_{{\mathbf x}_k} F \cdot \nabla_{{\mathbf p}_k} G -  \nabla_{{\mathbf x}_k} G \cdot \nabla_{{\mathbf p}_k} F \right) $$ 
and the following Hamiltonian system for each particle,
\begin{equation}
\label{eq:xpsym}
\begin{aligned}
& \dot{\mathbf x}_k = \frac{1}{w_k}\nabla_{{\mathbf p}_k}H, \\
& \dot{\mathbf p}_k =- \frac{1}{w_k}\nabla_{{\mathbf x}_k}H, \quad 1 \le k \le K,
\end{aligned}
\end{equation}
where $H$ is the discrete Hamiltonian~\eqref{eq:dish}.
Implicit symplectic mid-point rule is used to solve the above Hamiltonian system and preserve the symplectic structure~\cite{HLW}. Specifically, the scheme is 
\begin{equation}
\label{eq:xpmid}
\begin{aligned}
& \frac{{\mathbf x}_k^{n+1} - {\mathbf x}_k^{n}}{\Delta t} = {\mathbf p}_k^{n+\frac{1}{2}} - {\mathbf A}_h^{n}({\mathbf x}_k^{n+\frac{1}{2}}),\\
& \frac{{\mathbf p}_k^{n+1} - {\mathbf p}_k^{n}}{\Delta t} = - \left( \frac{\partial {\mathbf A}_h^{n}}{\partial \mathbf x}({\mathbf x}_k^{n+\frac{1}{2}}) \right)^\top  ({\mathbf A}^{n}_h({\mathbf x}_k^{n+\frac{1}{2}}) -  {\mathbf p}_k^{n+\frac{1}{2}}) + \text{term}_n,
\end{aligned}
\end{equation}
where 
\begin{equation}
\begin{aligned}
\text{term}_n & =   T  \sum_j w_j  \left(1 + \ln  \left(\sum_{k'=1}^K {w_{k'}} S({{\mathbf x} }_j - {{\mathbf x} }^{n+\frac{1}{2}}_{k'}) \right) \right)  \nabla S({{\mathbf x} }_j - {{\mathbf x} }^{n+\frac{1}{2}}_{k}).
\end{aligned}
\end{equation}
We use the Picard iteration to solve above implicit scheme, and denote the solution map of this subsystem as $\Phi_{xp}^{\Delta t}$.
\begin{remark}
As some derivatives are calculated in~\eqref{eq:xpmid} for the vector potential and smoothed delta functions, to guarantee the convergence of the iteration methods for solving the midpoint rule~\eqref{eq:xpmid}, degrees of B-splines in finite element space $V_0$ are at least $[3, 3, 3]$, and smoothed delta functions are second order B-splines at least. 
\end{remark}

In summary, we have the first and second order schemes for the first and second formulation,
\begin{align}
&\text{formulation I:} \quad \text{first order} \  \Phi_{xpa}^{\Delta t}\Phi_{a}^{\Delta t}, \ \text{second order} \  \Phi_{xpa}^{\Delta t/2}\Phi_{a}^{\Delta t} \Phi_{xpa}^{\Delta t/2},\label{eq:for1scheme} \\
&\text{formulation II:} \quad \text{first order} \ \Phi_{xp}^{\Delta t} \Phi_{a}^{\Delta t},   \ \text{second order} \ 
\Phi_{xp}^{\Delta t/2}\Phi_{a}^{\Delta t} \Phi_{xp}^{\Delta t/2}.\label{eq:for2scheme}
\end{align}
As the flow map of the formulation II~\eqref{eq:2ee2} does not preserve any term of the anti-symmetric bracket~\eqref{eq:original}, the schemes in~\eqref{eq:for2scheme} do not preserve it either.
However, schemes in~\eqref{eq:for2scheme} have better conservation properties and lower noise level than the schemes in~\eqref{eq:for1scheme}, as shown in the remaining parts of this paper.

\begin{remark}
For some simulations, there is a given background magnetic field ${\mathbf B}_0$. In this case, we should replace $\nabla \times {\mathbf A}$ with $\nabla \times {\mathbf A} + {\mathbf B}_0$ in the above subsystems, and the magnetic energy in Hamiltonian becomes $\frac{1}{2}\int |\nabla \times {\mathbf A} + {\mathbf B}_0|^2 \mathrm{d}{\mathbf x}$. Also we have another subsystem to solve, 
$$
 \frac{\partial f}{\partial t} = - ({\mathbf p} - {\mathbf A}) \times {\mathbf B}_0 \cdot \frac{\partial f}{\partial {\mathbf p}},
$$
which can be solved analytically. 
\end{remark}

\begin{remark}
The schemes constructed in this work can be applied to the case of adiabatic electrons, for which the $T\frac{\nabla n}{n}$ is replaced with $T\frac{\nabla n^\gamma}{n} = \frac{T \gamma}{\gamma-1}\nabla n^{\gamma-1}, \gamma \neq 1$.
\end{remark}

\begin{remark} 
Here we remark on other choices of the numerical methods of solving~\eqref{eq:xpsym}. In order to get energy conserving schemes, discrete gradient methods~\cite{DIS, go} can be used. 
Heavy Iterations are needed for solving the above implicit schemes, here we mention two explicit symplectic methods. 1) When the first order symplectic Euler method (explicit for ${\mathbf x}$ and implicit for ${\mathbf p}$) is applied for~\eqref{eq:xpsym}, it is in fact explicit as mentioned in~\cite{Qincanonical}, as we can update the ${\mathbf p}$ for each particle by inverting a $3 \times 3$ matrix and then update ${\mathbf x}$ explicitly. 2) High order explicit symplectic methods can be constructed by the explicit symplectic methods proposed in~\cite{ex1, ex2} combined with Hamiltonian splitting methods (for electron thermal energy term)~\cite{HLW, ha1}.
\end{remark}

\begin{remark}
When the electron effects is non-negligible as~\cite{valentini}, more complete Ohm's law should be used, i.e., 
$$
(I - d_e^2 \Delta ) {\mathbf E} = -({\mathbf u} \times {\mathbf B}) + \frac{1}{n} ({\mathbf J} \times {\mathbf B}) + \frac{1}{n} d_e^2 \nabla \cdot {\mathbf \Pi} - T \frac{\nabla n}{n} + \frac{d_e^2}{n} \nabla \cdot ({\mathbf u}{\mathbf J} + {\mathbf J}{\mathbf u} ) - \frac{1}{n} d_e^2 \nabla \cdot (\frac{{\mathbf{JJ}}}{n}),
$$
where ${\mathbf \Pi} = \int ({\mathbf v} - {\mathbf u}) ({\mathbf v} - {\mathbf u}) f \mathrm{d}{\mathbf v}$, and $d_e$ is the electron skin depth. We could also decompose the pressure term $T \frac{\nabla n}{n}$ as 
$$
T \frac{\nabla n}{n} = (I - d_e^2\Delta)\left(T \frac{\nabla n}{n}\right) + d_e^2 \Delta \left( T  \frac{\nabla n}{n}\right), 
$$
where the first term  $(I - d_e^2\Delta)T \frac{\nabla n}{n}$ gives the curl-free contribution of electric field ${\mathbf E}$. The same technique of moving $T \frac{\nabla n}{n}$ into Vlasov equation by changing the gauge used can also be used to reduce the effects of noise from particles. 
\end{remark}

\subsection{The comparisons of noise}
Here we comment on the noise of the above formulations I and II. As mentioned in sec.~\ref{sec:canonicalform}, in the formulation I~\eqref{eq:1ee}, the term $T \frac{\nabla n}{n}$ influences the update of the particles through the vector potential. Moreover, the derivative $ \frac{\partial {\mathbf A}}{\partial \mathbf x}$ in the Vlasov equation amplifies further the noise in ${\mathbf A}$. However, in the formulation II~\eqref{eq:2ee2}, the term $T \frac{\nabla n}{n}$ directly interacts with particles, and the noise in this term is not amplified as  in the formulation I. Then qualitatively, we know the schemes~\eqref{eq:for2scheme} of the formulation II~\eqref{eq:2ee2} are more accurate than the schemes~\eqref{eq:for1scheme} of the formulation I~\eqref{eq:1ee}. 
As the hybrid models are nonlinear, no explicit general solutions can be written out. In the test of finite grid instability with a special solution, an equilibrium,  in subsection~\ref{subsec:finite} , we compare the schemes of the two formulations quantitively. 

The schemes in~\eqref{eq:for1scheme} of the formulation I are improved when applying the filters introduced in the following for the density in the term $T \frac{\nabla n}{n}$, the results of which are shown in the numerical section~\ref{sec:numerical}. \\
\noindent{\bf Binomial filters}~\cite{BirdsallLangdon, Vay} The densities and the currents obtained from particles by deposition processes usually have large noise, a way to reduce the noise is to apply filters. The most commonly used filter in particle-in-cell simulations is the following three points filter
\begin{equation}\label{eq:filter}
\phi_j^f = \alpha  \phi_j + (1- \alpha) \frac{\phi_{j-1} + \phi_{j+1}}{2},
\end{equation}
where $\phi^f$ is the quantity after filtering. When $\alpha = 0.5$, it is called the binomial filter. 
When $\phi = e^{jkx}$, $\phi^f = g(\alpha, k) e^{jkx}$, where $g$ is called the filter gain, 
$$g(\alpha, k) = \alpha + (1-\alpha)\cos(k\Delta x) \approx 1 - (1-\alpha)\frac{(\Delta x)^2}{2} + \mathcal{O}(k^4).$$
When $m+1$ successive applications of filters with the coefficients $\alpha_1, \cdots, \alpha_m, \alpha_{m+1}$ are used, total attenuation $G$ is given by 
\begin{equation*}
G = \Pi_{i=1}^{m+1} g(\alpha_i, k) \approx 1 - \left( m + 1 - \sum_{i=1}^{m+1} \alpha_i\right) \frac{(k\Delta x)^2}{2} + \mathcal{O}(k^4).
\end{equation*}
The total attenuation $G$ is $1 + \mathcal{O}(k^4)$ if 
\begin{equation}\label{eq:filteralpha}
\alpha_{m+1} = m + 1 - \sum_{i=1}^{m}\alpha_i,
\end{equation}
where the $(m+1)$-th step is called a compensation step. Note that the values obtained after applying the binomial filters are non-negative if $\phi_j \ge 0$ for $\forall j$ in~\eqref{eq:filter}. The compensation step may produce negative values, although it does not happen in the numerical tests in next section. The application and investigation of non-negativity preserving filters deserve to be done in the future.

\section{Numerical experiments}\label{sec:numerical}
In this section, four numerical experiments are conducted to validate the codes of the above two schemes, and comparisons are made. For the scheme~\eqref{eq:for1scheme} of the formulation I~\eqref{eq:1ee}, it is illustrated that filters for the term $T\frac{\nabla n}{n}$ can improve the conservations of momentum and energy numerically. The scheme~\eqref{eq:for2scheme} for the formulation II~\eqref{eq:2ee2} is shown superior. The first order schemes in~\eqref{eq:for1scheme}-\eqref{eq:for2scheme} are used in the following, for both of which no filter is used for the subsystem~\eqref{eq:iseond}. The tolerance of Picard iteration
is set as $10^{-11}$. Periodic boundary conditions are considered. In the following, the 'with/without filtering' means if we use $m$ binomial filters with a compensation step satisfying~\eqref{eq:filteralpha} for the density $n$ of the term $T\frac{\nabla n}{n}$ in the scheme~\eqref{eq:for1scheme}, where $m$ is specified in each test. The implementations are conducted in the Python package STRUPHY~\cite{3}. 

\subsection{Finite grid instability}\label{subsec:finite}
\begin{figure}[htbp!]
\center{
\includegraphics[scale=0.4]{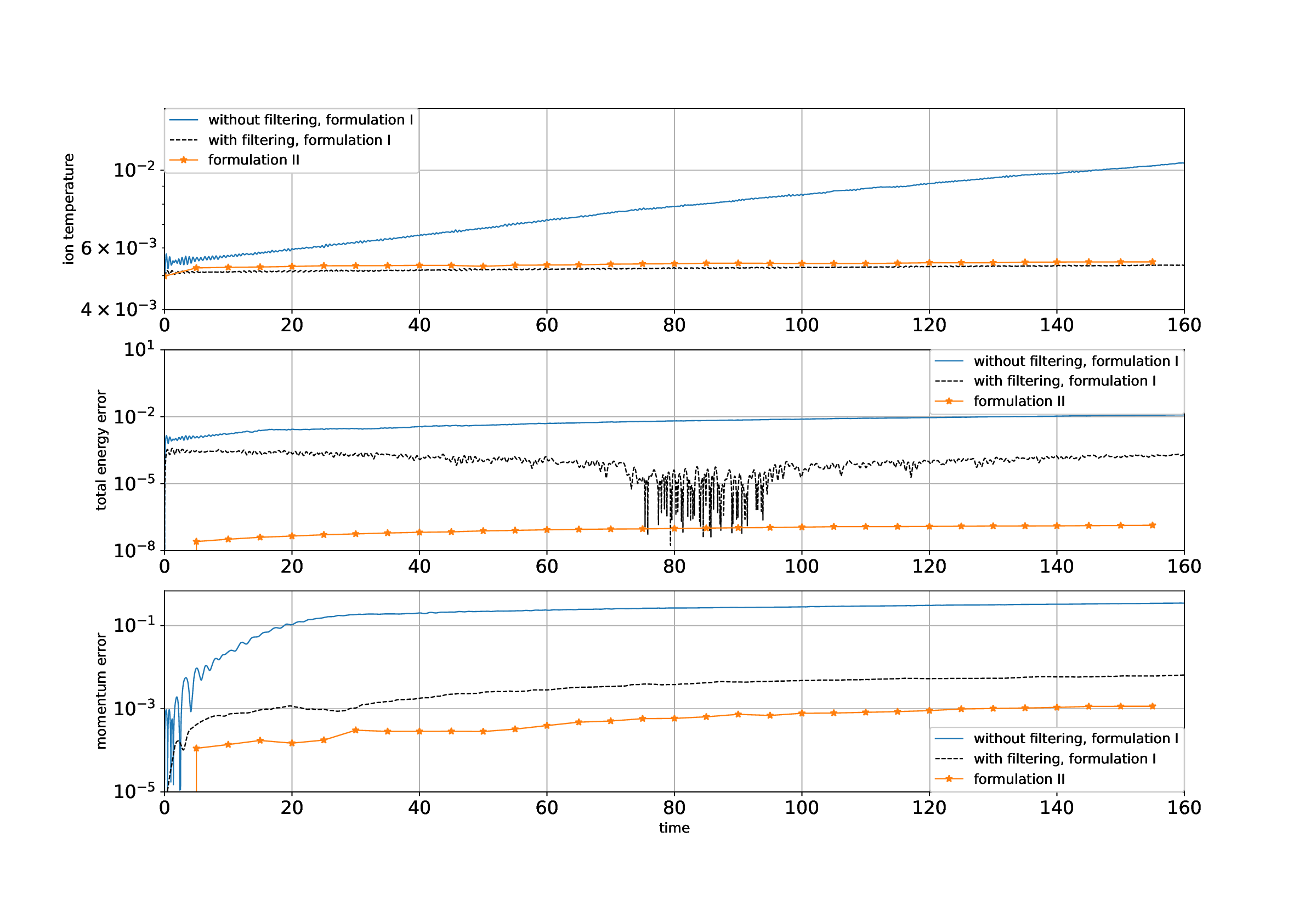}
}
\caption{{\bf Finite grid instability.}  Time evolutions of the ion temperature, relative energy error, and momentum error (the third component). }
 \label{finite_II}
\end{figure}

We validate our discretizations by a very challenging test called finite grid instability~\cite{Rambo}. In~\cite{chacon1}, the finite grid instability in the case of adiabatic electrons is studied numerically by a mass, momentum, energy conserving scheme, and the temperature of ions stay as a constant for a very long time.
In this test, a very cold ion beam with temperature $T_i = 0.005$ is moving with velocity $(0, 0, 0.1)$ in background electrons with temperature $T = 1$. Specifically, initial conditions are 
\begin{equation*}
\begin{aligned}
 {\mathbf B}_0 = {\mathbf 0}, \  {\mathbf A} = {\mathbf 0}, \ f =  \frac{1}{\pi^{\frac{3}{2}} v_T^{\frac{3}{2}}}  e^{-\frac{|p_x|^2}{ v_{T}^2} - \frac{|p_y|^2}{ v_{T}^2}  - \frac{|p_z - 0.1|^2}{ v_{T}^2} }, \,  T = 1, \ v_{T} = 0.1,
\end{aligned}
\end{equation*}
which is an equilibrium for the hybrid model, analytically it should stay unchanged with time. 
The computational domain is $[0, 1] \times [0, 1] \times [0, {5\pi}]$, the number of cells is $[4, 4, 32]$, degrees of B-splines are $[3,3,3]$, degrees of shape functions are $[2, 2, 2]$, quadrature points in each cell are $[2, 2, 4]$, total particle number is $5\times 10^4$, and the time step size is $0.01$. The first order schemes in~\eqref{eq:for1scheme} (without filter or with filters~\eqref{eq:filter} with $m = 8$ and coefficients satisfying~\eqref{eq:filteralpha}) and~\eqref{eq:for2scheme} are used. $m=8$ is large enough to filter out much noise in the density, and gives an almost constant ion temperature in Fig.~\ref{finite_II}. We tried to use filters with $m=4$ for scheme~\eqref{eq:for1scheme}, but we still can observe the obvious growth of the ion temperature. In this test for both schemes~\eqref{eq:for1scheme}-\eqref{eq:for2scheme} as $\nabla \times {\mathbf A}_h$ is 0,  $\Phi_{a}^{\Delta t}$ is just the identity map, and we are simply only using solution maps $\Phi_{xpa}^{\Delta t}$ and $\Phi_{xp}^{\Delta t}$ in schemes~\eqref{eq:for1scheme} and \eqref{eq:for2scheme}, respectively. 
From literature~\cite{Rambo}, we know that there is a quick growth of ion temperature with time when traditional particle-in-cell methods are used. 

\begin{figure}[htbp!]
\center{
\includegraphics[scale=0.4]{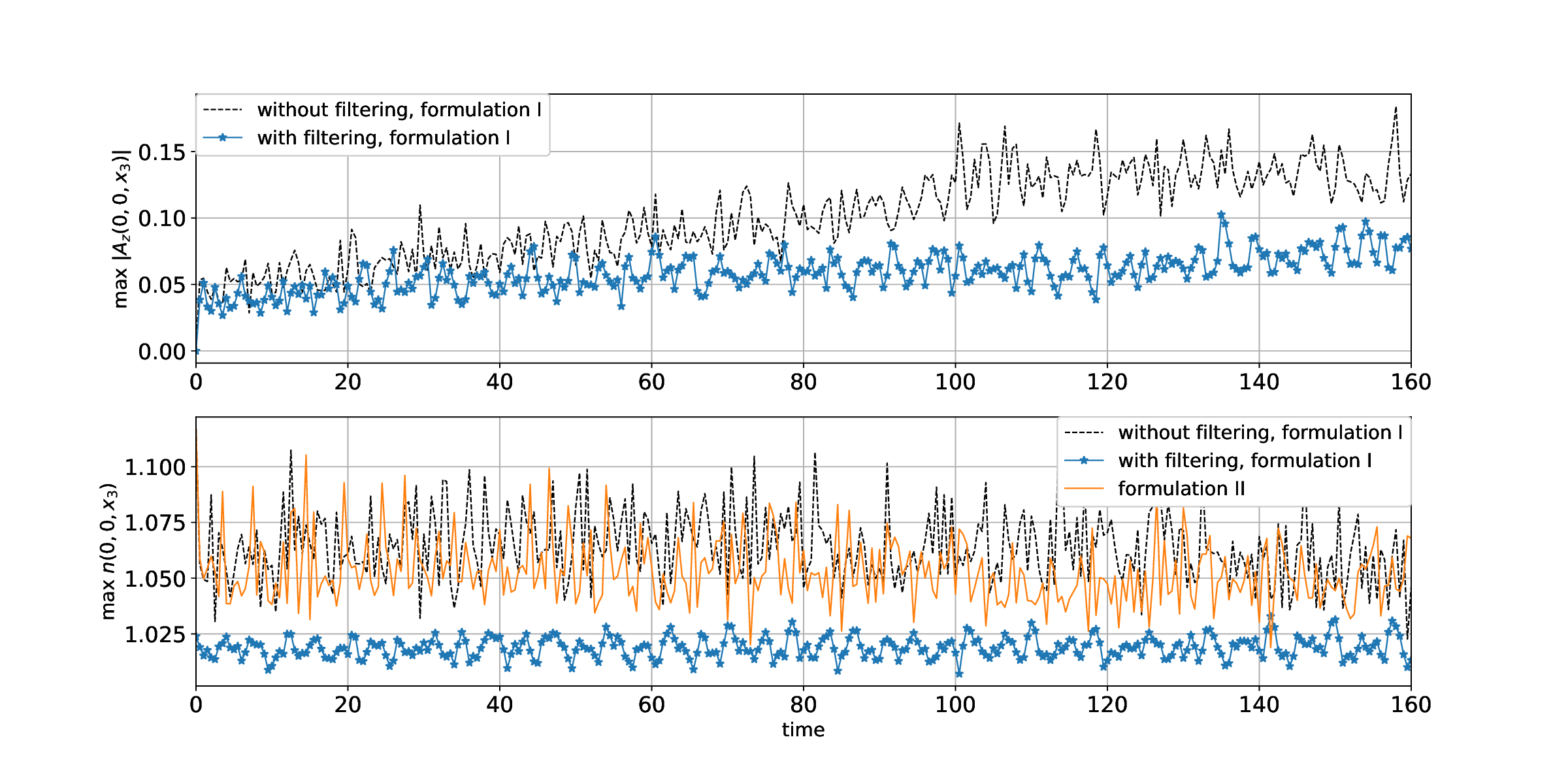}
}
\caption{{\bf Finite grid instability.} Time evolutions of the maximum values of $|A_z(0,0,x_3)|$, and density $n(0,0,x_3)$.}
 \label{finite_man}
\end{figure}

The following analysis about the noise holds when the solutions are close to the equilibrium. Let us firstly analyze the noise level of the formulation I~\eqref{eq:1ee}. The noise of the term $T \frac{\nabla n}{n}$ is at the level of $\mathcal{O}(E/\Delta s)$, where $E$ denotes the noise level of the density $n$ obtained from the particles and $\Delta s$ is the length of the cell in space.  By the the second equation in~\eqref{eq:ee}, we know the noise of the vector potential grows linearly with time, i.e., the noise of the vector potential would be $\mathcal{O}(tE/\Delta s)$ at time $t$.  The spatial derivatives of computing the term $ \frac{\partial {\mathbf A}}{\partial {\mathbf x}}$ amplify further the noises of the vector potential, and the noise of the term $\frac{\partial {\mathbf A}}{\partial \mathbf x}$ used in pushing particles is at the level of $\mathcal{O}(tE/(\Delta s)^2)$ at time $t$. 
However, in the formulation II, the vector potential is always zero, and the noise of the force of pushing particles, i.e., 
$T \frac{\nabla n}{n}$, is at the level of $\mathcal{O}(E/\Delta s)$. In Fig.~\ref{finite_man}, we plot the time evolutions of the maximum values of the third component of the vector potential ${\mathbf A}$ and the density. We can see that the maximum value of $|A_z(0,0,x_3)|$ given by the first order scheme~\eqref{eq:for1scheme} grows with time quickly when without filtering, and it grows much slower when filtering is used with a density closer to 1. The maximum value of $n(0,0,x_3)$ given by the scheme~\eqref{eq:for2scheme} of the formulation II oscillates around 1.05 with time. Note that the vector potential given by  the scheme~\eqref{eq:for2scheme} of the formulation II is always 0 in this test.

\begin{figure}[htbp!]
\center{
\includegraphics[scale=0.4]{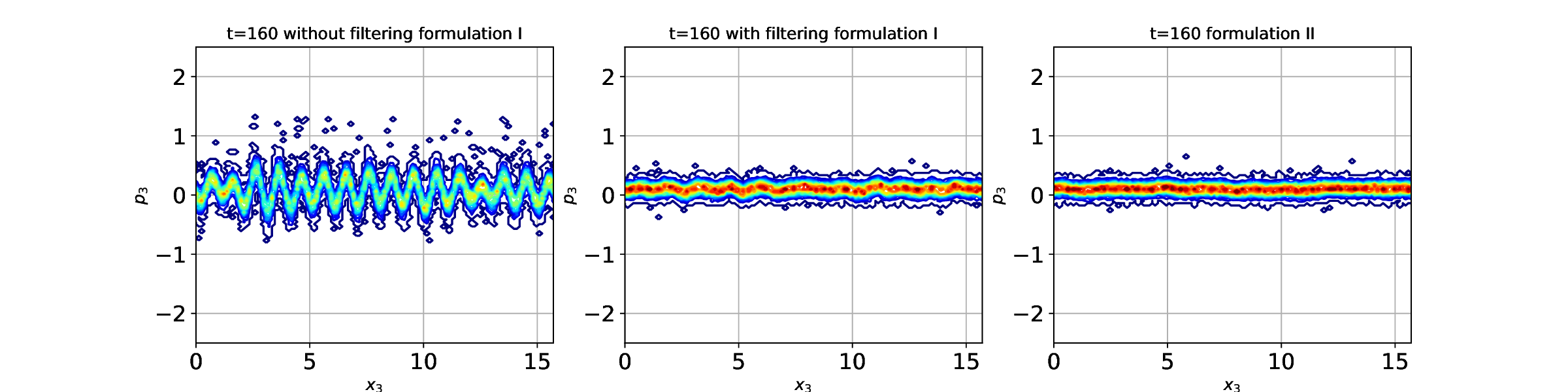}
}
\caption{{\bf Finite grid instability.} The contour plots of the scheme~\eqref{eq:for1scheme} with or without filtering and the scheme~\eqref{eq:for2scheme}.}
 \label{finite_I}
\end{figure}

From Fig.~\ref{finite_I}-\ref{finite_II}, we can see that first order scheme~\eqref{eq:for1scheme} without filtering gives an obvious growth of the temperature of ions with time, and the contour plot of $(x_3, p_3)$ at $t = 160$ is obviously distorted compared to the initial distribution function. However, when applying filters, the ion temperature grows much slower, and the contour plot at  $t = 160$ is still a very thin Maxwellian function.  Also we can see that the energy and momentum errors are much smaller when filters are used. The scheme~\eqref{eq:for2scheme} gives the best results, a thin Maxwellian at $t=160$, the smallest relative energy error around $10^{-7}$, and the smallest momentum error around $10^{-3}$. Note that during this simulation of the scheme~\eqref{eq:for2scheme}, no filter is used.


\subsection{Parallel electromagnetic wave: R mode} 
\begin{figure}[htbp!]
\center{
{\includegraphics[scale=0.4]{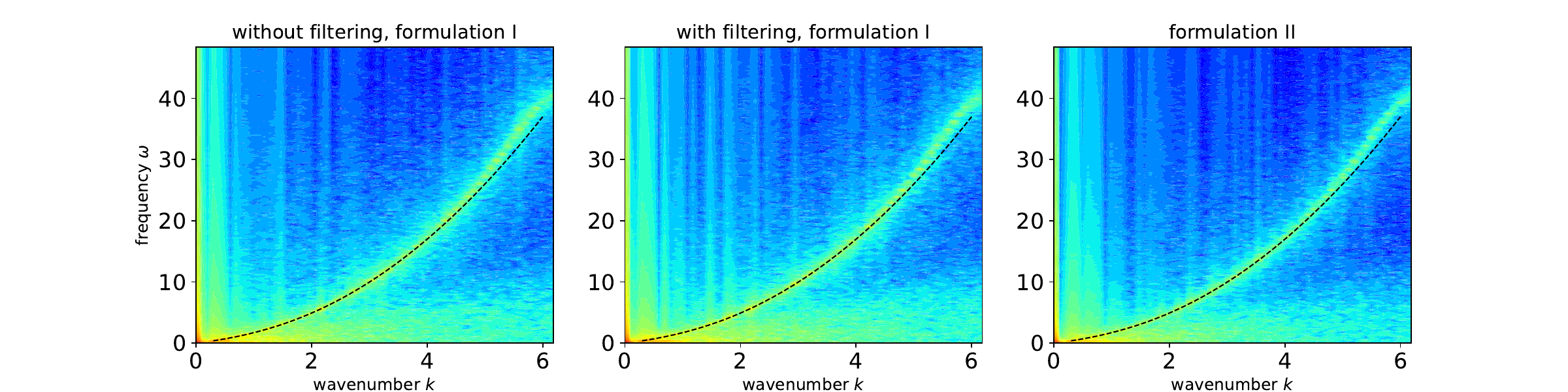}}
}
\caption{{\bf R-waves.}  The numerical dispersion relations of the R waves given by 1) the scheme~\eqref{eq:for1scheme} without filter, 2) the scheme~\eqref{eq:for1scheme} with filters, 3) the scheme~\eqref{eq:for2scheme}.}
 \label{R4}
\end{figure}

\begin{figure}[htbp!]
\center{
{\includegraphics[scale=0.4]{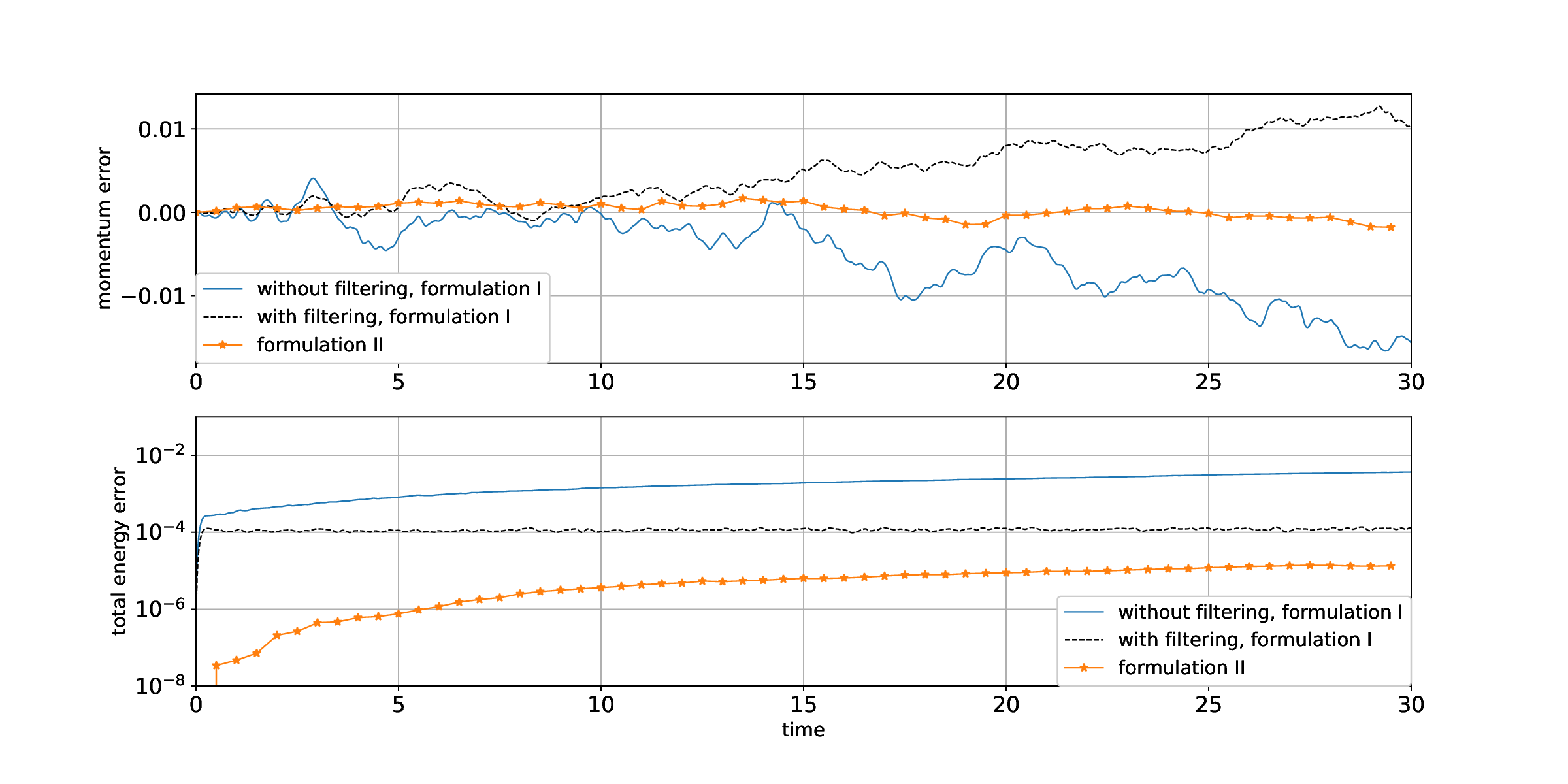}}
}
\caption{{\bf R-waves.} Time evolutions of the momentum errors (the third component) and relative energy errors given by 1) the scheme~\eqref{eq:for1scheme} without filter, 2) the scheme~\eqref{eq:for1scheme} with filters, 3) the scheme~\eqref{eq:for2scheme}.}
 \label{R1}
\end{figure}

Then we check a parallel propagating R wave by a quasi-1D simulation. Background magnetic field is along $z$ direction. No perturbation is added for the system other than the noise of PIC due to the reduced number of macro-particles. Specifically, the initial conditions we use are:
\begin{equation*}
\begin{aligned}
{\Bb}_0 = (0, 0, 1), \, {\mathbf A} = (0, 0, 0), \, T  = 1, \, f =  \frac{1}{\pi^{\frac{3}{2}} }e^{-{|{\mathbf p}|^2} }.
\end{aligned}
\end{equation*}
Computational parameters are: 
grid number $[4, 4, 128]$, domain $[0, 1] \times [0, 1] \times [0, 64]$, $dt = 0.005$, final computation time $40$, total particle number $2\times 10^{5}$, degree of B-splines $[3, 3, 3]$, quadrature point in each cell $[2,2,4]$, and degrees of shape functions $[2,2,2]$. The cell sizes in three directions are $[0.25, 0.25, 0.5]$, which are small enough to prevent the large errors due to the cancellation problem as pointed in~\cite{cancellation} for high-order shape functions. The first order schemes in~\eqref{eq:for1scheme} (without filter or with filters~\eqref{eq:filter} with $m=2$ and coefficients satisfying~\eqref{eq:filteralpha}) and~\eqref{eq:for2scheme} are used.
See the numerical results of dispersion relation of R mode given by the scheme~\eqref{eq:for1scheme} in Fig.~\ref{R4}. The black dash lines are the analytical dispersion relations given by Python package HYDRO proposed in \cite{Told}, when $k \ll 1$, $\omega \propto k$, when $k \gg 1$, $\omega \propto k^2$. We can see that our numerical results of the schemes~\eqref{eq:for1scheme}-\eqref{eq:for2scheme} are in good agreements with the analytical results pretty well even when the wave number $k$ is larger than Nyquist frequency. In Fig.~\ref{R1}, we present the time evolutions of momentum errors and relative energy errors of the schemes~\eqref{eq:for1scheme}-\eqref{eq:for2scheme}, from which we can see the scheme~\eqref{eq:for2scheme} is better than the scheme~\eqref{eq:for1scheme}. For the scheme~\eqref{eq:for1scheme}, the relative energy error with filters is smaller than the results without filters. The momentum error and relative energy error of the scheme~\eqref{eq:for2scheme} without filter are the smallest, the relative energy error is at the level of $10^{-6}$.

\subsection{Perpendicular wave: ion Bernstein waves}
\begin{figure}[htbp!]
\center{
{\includegraphics[scale=0.4]{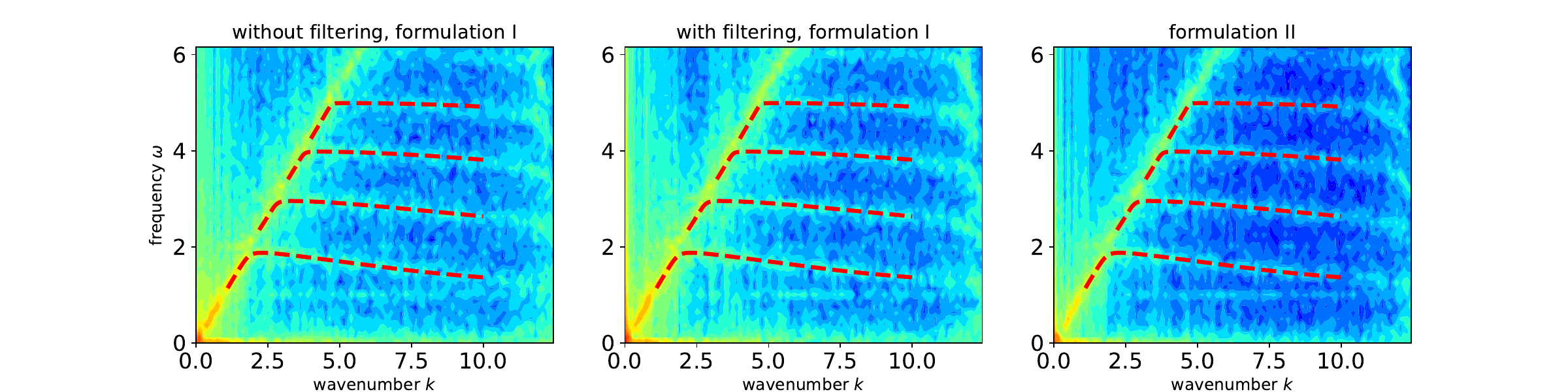}}
}
\caption{{\bf Bernstein waves.} Bernstein waves dispersion relation of (a) the scheme~\eqref{eq:for1scheme} without filtering; (b) the scheme~\eqref{eq:for1scheme} with filtering; 3) the scheme~\eqref{eq:for2scheme}. Red dash lines are different branches of the analytical dispersion relation of Bernstein waves.}
 \label{B1}
\end{figure}

\begin{figure}[htbp!]
\center{
{\includegraphics[scale=0.4]{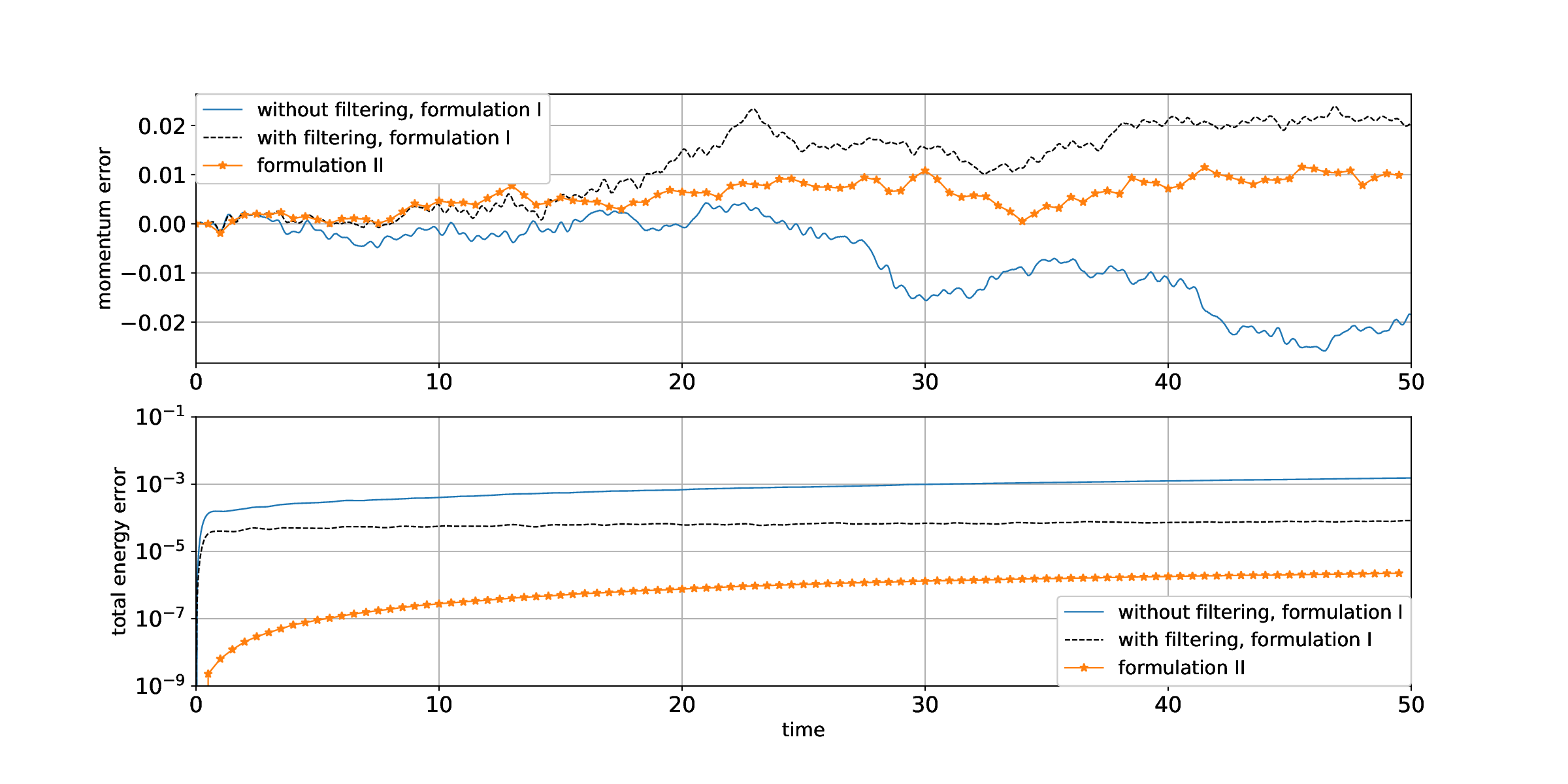}}
}
\caption{{\bf Bernstein waves.} Time evolutions of the relative energy errors and momentum errors (the first component) given by the scheme~\eqref{eq:for1scheme} and scheme~\eqref{eq:for2scheme}.}
 \label{B2}
\end{figure}

Then we check Bernstein waves by a one dimensional simulation, which are perpendicular to background magnetic field. 
In order to excite these waves, we initialize a quasi-1D thermal plasma along the $x$ direction. No initial perturbation is added except the noise of PIC method. Specifically, initial conditions are:
\begin{equation*}
\begin{aligned}
 {\Bb}_0 = (0, 0, 1), \,  {\mathbf A} = (0, 0, 0), \ f =  \frac{1}{\pi^{\frac{3}{2}} v_T^{\frac{3}{2}}} e^{-\frac{|{\mathbf p}|^2}{v_T^2} }, \, \green{T}  = 0.09. 
\end{aligned}
\end{equation*}
Computational parameters are: 
grid number $[200, 4, 4]$, domain $[0, 50] \times [0, 1] \times [0, 1]$, time step size $0.005$, $v_T = 0.2121$, final computation time $80$, particle number $10^5$, degrees of polynomials $[3,3,3]$, quadrature point in each cell $[4,2,2]$, and degrees of shape functions $[2,2,2]$. The first order schemes in~\eqref{eq:for1scheme} (without filter or with filters with $m=3$ and  coefficients satisfying~\eqref{eq:filteralpha}) and~\eqref{eq:for2scheme} are used.
Firstly, we check the numerical dispersion relations given by the scheme~\eqref{eq:for1scheme} and scheme~\eqref{eq:for2scheme}, which are presented in Fig.~\ref{B1}. We can see that both schemes give the correct dispersion relation numerically,  the red dashed lines are analytical dispersion relations of Bernstein waves obtained via HYDRO code~\cite{Told}.  
All the $n$-th band of the dispersion relation starts close to $n+1$ for small wavenumber and end up at $n$ for large wavenumber. In Fig.~\ref{B2}, we present the momentum errors and relative energy errors of the two schemes.  For the scheme~\eqref{eq:for1scheme},  the relative energy error and momentum error become smaller when filters are used.  The scheme~\eqref{eq:for2scheme} gives best results, the relative energy error is conserved at the level of $10^{-6}$, and the momentum error is around $10^{-2}$. 

\subsection{Ion cyclotron instability}
Finally we validate the numerical schemes via a nonlinear instability test, ion cyclotron instability. 
We initialize a quasi-1D thermal plasma along the $x$ direction. No initial perturbation is added except the noise of PIC method. Specifically, initial conditions are:
\begin{equation*}
\begin{aligned}
 {\Bb}_0 = (1, 0, 0), \,  {\mathbf A} = (0, 0, 0), \ f =  \frac{1}{\pi^{\frac{3}{2}} v_\parallel v_\perp^2} e^{-\frac{p_x^2}{v_\parallel^2} -\frac{p_y^2+p_z^2}{v_\perp^2} }, \, T  = 0.5,
\end{aligned}
\end{equation*}
in which there is a temperature anisotropy. 
Computational parameters are: 
grid number $[32, 4, 4]$, domain $[0, \frac{2\pi}{0.7}] \times [0, 1] \times [0, 1]$, time step size $0.005$, $v_\parallel = 1$, $v_\perp = 2$, final computation time $100$, degrees of polynomials $[3,3,3]$, quadrature point in each cell $[4,2,2]$, and degrees of shape functions $[2,2,2]$. 
The first order schemes in~\eqref{eq:for1scheme} (without filter or with filters with $m=1$ and coefficients satisfying~\eqref{eq:filteralpha}) and~\eqref{eq:for2scheme} are used.

\begin{figure}[htbp!]
\center{
{\includegraphics[scale=0.42]{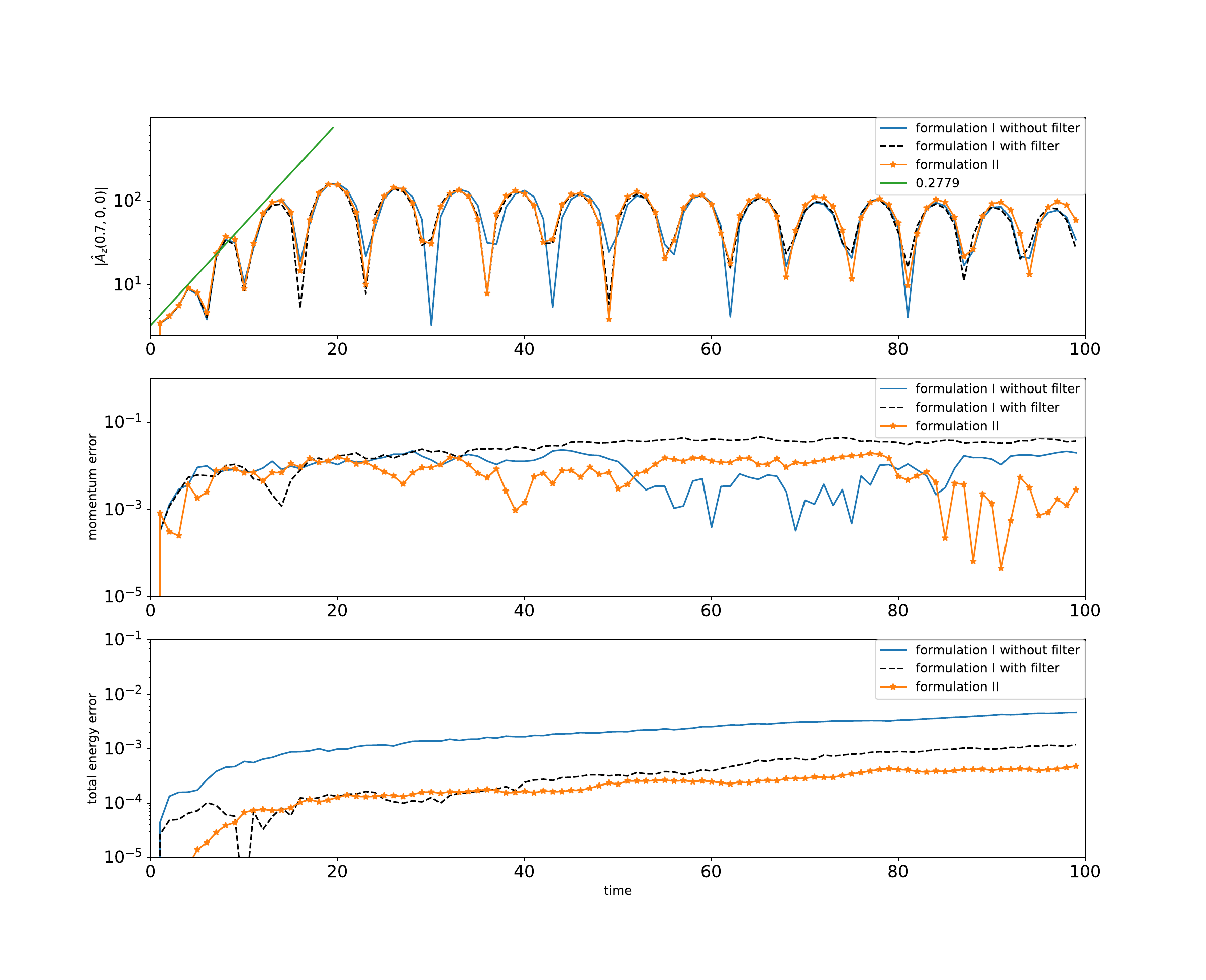}}
}
\caption{{\bf Ion cyclotron instability.} Time evolutions of the amplitude of the Fourier mode $\hat{A}_z(0.7, 0, 0)$, momentum error (the first component), and relative energy error with $2\times 10^4$ particles.}
 \label{ion5000}
\end{figure}

We present the numerical results of the schemes~\eqref{eq:for1scheme}-\eqref{eq:for2scheme} using 20000 particles in total in Fig.~\ref{ion5000}. We can see both schemes give correct exponential growth rates of the Fourier model $\hat{A}_z(t, 0.7, 0, 0)$, 0.2779, which is obtained by using the dispersion relation solver HYDRO code~\cite{Told}. For the scheme~\eqref{eq:for1scheme} without the filter, the crest of the last oscillations of $|\hat{A}_z(t, 0.7, 0, 0)|$ at about $t = 100$ is smallest compared to the results of the scheme~\eqref{eq:for1scheme} with filters and the scheme~\eqref{eq:for2scheme}. The scheme~\eqref{eq:for2scheme} gives the smallest momentum error (when time is around 90) and relative energy error. 
The scheme~\eqref{eq:for1scheme} gives smaller energy error when filters are used. 

\section{Conclusion}\label{sec:conclusion}

In this work, we explore and compare the canonical momentum based particle-in-cell methods for the two formulations  of the hybrid model with kinetic ions and massless electrons. Splitting methods and mid-point rules are used for time discretizations. The schemes of the first formulation are significantly improved by using binomial filters. The schemes of the second formulation show better conservation properties even without filter. Finding more efficient solvers for~\eqref{eq:picardaa}, applying our numerical methods to large scale physical simulations, and numerical comparisons with the schemes constructed in~\cite{LHPS} are future works. Smoothed delta functions are used only for getting the discrete density on the grids and calculating the electron thermal energy, which can also be used in the discretizations of the other particle related parts. 
Let us also mention that other than finite element methods, other methods, such as spectral methods~\cite{spectral} can be used to discretize the vector potential. Decomposing the right hand side of the Faraday's law in other forms, such as using Hodge-Helmholtz decompostion and other gauges deserves, to be investigated.

\section*{Acknowledgements}
Simulations in this work were performed on Max Planck Computing \& Data Facility (MPCDF). 
The authors would like to thank anonymous reviewers for many helpful comments for improving this paper.

\section{Appendix}
\subsection{The hybrid model with a background magnetic field}\label{subsec:back}
When there is a background magnetic field ${\mathbf B}_0$, the formulation I~\eqref{eq:1ee} is
\begin{equation*}
\begin{aligned}
& \frac{\partial f}{\partial t} = -({\mathbf p} - {\mathbf A}) \cdot \frac{\partial f}{\partial {\mathbf x}} + \left[\left( \frac{\partial {\mathbf A}}{\partial \mathbf x} \right)^\top  ({\mathbf A} - {\mathbf p})\right] \cdot \frac{\partial f}{\partial {\mathbf p}} - ({\mathbf p} - {\mathbf A}) \times {\mathbf B}_0 \cdot \frac{\partial f}{\partial {\mathbf p}},\\
& \frac{\partial {\mathbf A}}{\partial t} = T \frac{\nabla n}{n} - \frac{\nabla \times \left({\mathbf B}_0 +  \nabla \times {\mathbf A}\right)}{n} \times \left( {\mathbf B}_0 + {\nabla \times {\mathbf A}}\right) - \frac{\int ({\mathbf A} - {\mathbf p}) f \mathrm{d}{\mathbf p}}{n} \times \left({{\mathbf B}_0 + \nabla \times {\mathbf A}}\right).
\end{aligned}
\end{equation*}
The formulation II~\eqref{eq:2ee2} becomes
\begin{equation*}
\begin{aligned}
& \frac{\partial f}{\partial t} = -({\mathbf p} - {\mathbf A}) \cdot \frac{\partial f}{\partial {\mathbf x}} + \left[T \frac{\nabla n}{n}  + \left( \frac{\partial {\mathbf A}}{\partial \mathbf x} \right)^\top  ({\mathbf A} - {\mathbf p})\right] \cdot \frac{\partial f}{\partial {\mathbf p}} - ({\mathbf p} - {\mathbf A}) \times {\mathbf B}_0 \cdot \frac{\partial f}{\partial {\mathbf p}},\\
& \frac{\partial {\mathbf A}}{\partial t} = - \frac{\nabla \times \left({\mathbf B}_0 +  \nabla \times {\mathbf A}\right)}{n} \times \left( {\mathbf B}_0 + {\nabla \times {\mathbf A}}\right) - \frac{\int ({\mathbf A} - {\mathbf p}) f \mathrm{d}{\mathbf p}}{n} \times \left({{\mathbf B}_0 + \nabla \times {\mathbf A}}\right),
\end{aligned}
\end{equation*}
which could be derived with the following Hamiltonian and anti-symmetric bracket,
\begin{equation*}
\mathcal{H} = \frac{1}{2}   \int f |{\mathbf p} - {\mathbf A}|^2   \mathrm{d}{\mathbf x} \mathrm{d}{\mathbf p} + T \int n \ln n \mathrm{d}{\mathbf x} + \frac{1}{2} \int |\nabla \times {\mathbf A} + {\mathbf B}_0|^2   \mathrm{d}{\mathbf x},
\end{equation*}
\begin{equation*}
\begin{aligned}
\{ \mathcal{F}, \mathcal{G} \}(f, {\mathbf A}) &=   \int f   \left[\frac{\delta \mathcal{F}}{\delta f}, \frac{\delta \mathcal{G}}{\delta f} \right]_{xp}   \mathrm{d}{\mathbf x} \mathrm{d}{\mathbf p}  - \int \frac{1}{n} \left(\nabla \times \mathbf{A}+ {\mathbf B}_0\right) \cdot \left(  \frac{\delta \mathcal{F}}{\delta {\mathbf A}} \times \frac{\delta \mathcal{G}}{\delta {\mathbf A}} \right) \mathrm{d}{\mathbf x} \\
& + \int f {\mathbf B}_0 \cdot \left( \frac{\partial }{\partial {\mathbf p}} \frac{\delta {\mathcal{F}}}{\delta f} \times  \frac{\partial }{\partial {\mathbf p}} \frac{\delta {\mathcal{G}}}{\delta f}\right) \mathrm{d}{\mathbf p} \mathrm{d}{\mathbf x}.
\end{aligned}
\end{equation*}

\subsection{$xvA$ formulation}\label{eq:xvAbracket}
The equivalent hybrid model with unknowns $f(t, {\mathbf x}. {\mathbf v})$ and ${\mathbf A}(t,{\mathbf x})$ is 
\begin{equation*}
\begin{aligned}
& \frac{\partial f}{\partial t} + {\mathbf v} \cdot \frac{\partial f}{\partial \xb} +({\mathbf E} + {\mathbf v} \times {\Bb}) \cdot \frac{\partial f}{\partial {\mathbf v}} = 0\,,\\
&\frac{\partial {\mathbf A}}{\partial t} = \left(\mathbf u - \frac{\mathbf J}{n}  \right) \times {\Bb}\,,\\
& n = \int f \mathrm{d}{\mathbf v}, \quad n{\mathbf u} = \int {\mathbf v}f \mathrm{d}{\mathbf v}, \quad \Bb = \nabla \times {\mathbf A}\,,\\
& {\mathbf E} = - T \frac{\nabla n}{n} - \left(\mathbf u - \frac{\mathbf J}{n}  \right) \times {\Bb}, \quad \mathbf J = \nabla \times {\mathbf B}.\\
\end{aligned}
\end{equation*}
We propose the following anti-symmetric bracket and total energy for this formulation.
\begin{equation*}
\begin{aligned}
\{ \mathcal{F}, \mathcal{G} \}(f, {\mathbf A}) &= \int f  \left[\frac{\delta \mathcal{F}}{\delta f}, \frac{\delta \mathcal{G}}{\delta f} \right]_{xv}   \mathrm{d}{\mathbf x} \mathrm{d}{\mathbf v} + \int \frac{\nabla \times {\mathbf A}}{\int f \mathrm{d}{\mathbf v}} \cdot \frac{\delta \mathcal{F}}{\delta {\mathbf A}} \times \frac{\delta \mathcal{G}}{\delta {\mathbf A}} \mathrm{d} {\mathbf x}\\
& + \int f \frac{\nabla \times {\mathbf A}}{\int f \mathrm{d}{\mathbf v}} \cdot \left( \frac{\partial }{\partial {\mathbf v}} \frac{\delta \mathcal{G}}{\delta f} \times \frac{\delta \mathcal{F}}{\delta {\mathbf A}} -  \frac{\partial }{\partial {\mathbf v}} \frac{\delta \mathcal{F}}{\delta f} \times \frac{\delta \mathcal{G}}{\delta {\mathbf A}}   \right) \mathrm{d}{\mathbf x} \mathrm{d}{\mathbf v}\\
& + \int \frac{\nabla \times {\mathbf A}}{\int f \mathrm{d}{\mathbf v}} \cdot \left( \int f \frac{\partial }{\partial {\mathbf v}} \frac{\delta \mathcal{F}}{\delta f}  \mathrm{d}{\mathbf v}\right) \times \left( \int f \frac{\partial }{\partial {\mathbf v}'} \frac{\delta \mathcal{G}}{\delta f}  \mathrm{d}{\mathbf v}'\right)  \mathrm{d}{\mathbf x} \\
& + \int f \left( \nabla \times {\mathbf A} \cdot \frac{\partial}{\partial {\mathbf v}} \frac{\delta \mathcal{F}}{\delta f} \times  \frac{\partial}{\partial {\mathbf v}} \frac{\delta \mathcal{G}}{\delta f}  \right) \mathrm{d}{\mathbf x} \mathrm{d}{\mathbf v},
\end{aligned}
\end{equation*}
$$
\cH =  \frac{1}{2} \int |{\mathbf v}|^2f\, \tn{d}{\xb} \tn{d}{\mathbf v}
 + \frac{1}{2} \int |\nabla \times {\mathbf A}|^2\, \tn{d}{\xb}\, +  T \int \left( \int f\,\tn d \vb\right) \ln \left( \int f\,\tn d \vb\right) \tn{d}{\xb}.
 $$
Energy conserving particle-in-cell methods for the $xvA$ formulation can be constructed following the strategies proposed in~\cite{2020en, LHPS}.


\begin{thebibliography}{00}

\bibitem{GEMPIC} Kraus M, Kormann K, Morrison P J, Sonnendr\"ucker E. GEMPIC: geometric electromagnetic particle-in-cell methods. Journal of Plasma Physics, 2017, 83(4).

\bibitem{Qincanonical} Qin H, Liu J, Xiao J, Zhang R, He Y, Wang Y, Sun Y, Burby J W, Ellison L, Zhou Y. Canonical symplectic particle-in-cell method for long-term large-scale simulations of the Vlasov--Maxwell equations. Nuclear Fusion, 2015, 56(1): 014001.

\bibitem{1}{Tronci C. Hamiltonian approach to hybrid plasma models. Journal of Physics A, 2010, 43(37).}

\bibitem{DECVM} Xiao J, Qin H, Liu J, He Y, Zhang R, Sun Y. Explicit high-order non-canonical symplectic particle-in-cell algorithms for Vlasov-Maxwell systems. Physics of Plasmas, 2015, 22(11): 112504.


\bibitem{Feng} Feng K, Qin M. Symplectic geometric algorithms for Hamiltonian systems. Berlin: Springer, 2010.

\bibitem{3} {Holderied F, Possanner S, Wang X. MHD-kinetic hybrid code based on structure-preserving finite elements with particles-in-cell. Journal of Computational Physics, 2021, 433: 110143.}

\bibitem{qinxiao} Xiao J, Qin H. Field theory and a structure-preserving geometric particle-in-cell algorithm for drift wave instability and turbulence. Nuclear Fusion, 2019, 59(10): 106044.

\bibitem{BirdsallLangdon} Birdsall C K, Langdon A B. Plasma physics via computer simulation. CRC press, 2018.

 
\bibitem{HLW} Hairer E, Lubich C, Wanner G. Geometric Numerical Integration: Structure-Preserving Algorithms for Ordinary Differential Equations, vol. 31, Springer Science \& Business Media, 2006.

\bibitem{FEEC} Arnold D N, Falk R S, Winther R. Finite element exterior calculus, homological techniques, and applications. Acta numerica, 2006, 15: 1-155.

\bibitem{DEC} Hirani A N. Discrete exterior calculus. California Institute of Technology, 2003.


\bibitem{2020en}Kormann K, Sonnendr\"ucker E. Energy-conserving time propagation for a structure-preserving particle-in-cell Vlasov--Maxwell solver. Journal of Computational Physics, 2020, 425: 109890.

\bibitem{hevm} He Y, Sun Y, Qin H, Liu J. Hamiltonian particle-in-cell methods for Vlasov--Maxwell equations. Physics of Plasmas, 2016, 23(9): 092108.


\bibitem{iter} Kelley C T. Iterative methods for linear and nonlinear equations. Society for Industrial and Applied Mathematics, 1995.

\bibitem{go} Gonzalez O. Time integration and discrete Hamiltonian systems. Journal of Nonlinear Science, 1996, 6(5): 449-467.
 
\bibitem{kaltsas} Kaltsas D A, Throumoulopoulos G N, Morrison P J. Hamiltonian kinetic-Hall magnetohydrodynamics with fluid and kinetic ions in the current and pressure coupling schemes. Journal of Plasma Physics, 2021, 87(5): 835870502.
 
\bibitem{Matthews} Matthews A P. Current advance method and cyclic leapfrog for 2D multispecies hybrid plasma simulations. Journal of Computational Physics, 1994, 112(1): 102-116.
 
\bibitem{CAMELIA} Franci L, Hellinger P, Guarrasi M, Chen C H K,  Papini E, Verdini A, Matteini L, Landi S. Three-dimensional simulations of solar wind turbulence with the hybrid code CAMELIA. Journal of Physics: Conference Series. IOP Publishing, 2018, 1031(1): 012002.

\bibitem{valentini} Valentini F, Tr\'{a}vn\'{i}\v{c}ek P, Califano F, Hellinger P, Mangeney A. A hybrid-Vlasov model based on the current advance method for the simulation of collisionless magnetized plasma. Journal of Computational Physics, 2007, 225(1): 753-770.


\bibitem{LHPS} Li Y, Campos Pinto M, Holderied F, Possanner S, Sonnendr\"ucker E.  Geometric Particle-In-Cell discretizations of a plasma hybrid model with kinetic ions and mass-less fluid electrons. Journal of Computational Physics, 2024, 498: 112671.

\bibitem{chacon1} Stanier A, Chac\'on L, Chen G. A fully implicit, conservative, non-linear, electromagnetic hybrid particle-ion/fluid-electron algorithm. Journal of Computational Physics, 2019, 376: 597-616.

\bibitem{chacon2} Stanier A, Chac\'on L. A conservative implicit-PIC scheme for the hybrid kinetic-ion fluid-electron plasma model on curvilinear meshes. Journal of Computational Physics, 2022, 459: 111144.

\bibitem{Told} Told D, Cookmeyer J, Astfalk P, Jenko F. A linear dispersion relation for the hybrid kinetic-ion/fluid-electron model of plasma physics. New Journal of Physics, 2016, 18(7): 075001.

\bibitem{Pegasus} Kunz M W, Stone J M, Bai X N. Pegasus: a new hybrid-kinetic particle-in-cell code for astrophysical plasma dynamics. Journal of Computational Physics, 2014, 259: 154-174.

\bibitem{Vay} Vay J L, Geddes C G R, Cormier-Michel E, Grote D P. Numerical methods for instability mitigation in the modeling of laser wakefield accelerators in a Lorentz-boosted frame. Journal of Computational Physics, 2011, 230(15): 5908-5929.

\bibitem{Rambo} Rambo P W. Finite-grid instability in quasineutral hybrid simulations. Journal of Computational Physics, 1995, 118(1): 152-158.

\bibitem{23} Lipatov A S. The hybrid multiscale simulation technology: an introduction with application to astrophysical and laboratory plasmas. Springer Science \& Business Media, 2002.

\bibitem{26} Park W, Parker S, Biglari H, et al. Three-dimensional hybrid gyrokinetic-magnetohydrodynamics simulation. Physics of Fluids B: Plasma Physics, 1992, 4(7): 2033-2037.

\bibitem{31} Winske D, Yin L, Omidi N, Karimabadi H, Quest K. Hybrid simulation codes: Past, present and future--A tutorial. Space plasma simulation, 2003: 136-165.


\bibitem{DIS} McLachlan R I, Quispel G R W, Robidoux N. Geometric integration using discrete gradients. Philosophical Transactions of the Royal Society of London. Series A: Mathematical, Physical and Engineering Sciences, 1999, 357(1754): 1021-1045.

\bibitem{ex1} Zhang R, Qin H, Tang Y, Liu J, He Y, Xiao J. Explicit symplectic algorithms based on generating functions for charged particle dynamics. Physical Review E, 2016, 94(1): 013205.

\bibitem{ex2} Zhou Z, He Y, Sun Y, Liu J, Qin H. Explicit symplectic methods for solving charged particle trajectories. Physics of Plasmas, 2017, 24(5): 052507.

\bibitem{rehy} Haggerty C C, Caprioli D. dHybridR: A hybrid particle-in-cell code including relativistic ion dynamics. The Astrophysical Journal, 2019, 887(2): 165.

\bibitem{frame} Campos Pinto M, Kormann K, Sonnendr\"ucker E. Variational framework for structure-preserving electromagnetic particle-in-cell methods. Journal of Scientific Computing, 2022, 91(2): 1-39.


\bibitem{Winske} Winske D. Hybrid simulation codes with application to shocks and upstream waves. Space Science Reviews, 1985, 42(1): 53-66.

\bibitem{newreview} Winske D, Karimabadi H, Le A, Omidi N, Roytershteyn V, Stanier A. Hybrid codes (massless electron fluid). arXiv preprint arXiv:2204.01676, 2022.

\bibitem{chenliu} Chen L, White R B, Rosenbluth M N. Excitation of internal kink modes by trapped energetic beam ions. Physical Review Letters, 1984, 52(13): 1122.

\bibitem{Porcelli} Coppi B, Porcelli F. Theoretical model of fishbone oscillations in magnetically confined plasmas. Physical review letters, 1986, 57(18): 2272.

\bibitem{guoyongfu} Park W, Belova E V, Fu G Y, Tang X Z, Strauss H R, Sugiyama L E. Plasma simulation studies using multilevel physics models. Physics of Plasmas, 1999, 6(5): 1796-1803.

\bibitem{pham}  Hahm T S, Lee W W, Brizard A. Nonlinear gyrokinetic theory for finite-beta plasmas. The Physics of fluids, 1988, 31(7): 1940-1948.

\bibitem{Mishchenko} Mishchenko A, Cole M, Kleiber R, K\"onies A. New variables for gyrokinetic electromagnetic simulations. Physics of Plasmas, 2014, 21(5).

\bibitem{dongjian} Bao J, Liu D, Lin Z. A conservative scheme of drift kinetic electrons for gyrokinetic simulation of kinetic-MHD processes in toroidal plasmas. Physics of Plasmas, 2017, 24(10).

\bibitem{LinZ} Bao J, Lin Z, Lu Z X. A conservative scheme for electromagnetic simulation of magnetized plasmas with kinetic electrons. Physics of Plasmas, 2018, 25(2).

\bibitem{spectral} Campos Pinto M, Ameres J, Kormann K, Sonnendr\"ucker E. On geometric Fourier particle in cell methods. arXiv preprint arXiv:2102.02106, 2021.

\bibitem{valentine} Servidio S, Valentini F, Califano F, Veltri P. Local kinetic effects in two-dimensional plasma turbulence. Physical review letters, 2012, 108(4): 045001.

\bibitem{Schekochihin} Kunz M W, Schekochihin A A, Stone J M. Firehose and mirror instabilities in a collisionless shearing plasma. Physical Review Letters, 2014, 112(20): 205003.

\bibitem{Arzamasskiy} Arzamasskiy L, Kunz M W, Chandran B D G, Quataert E. Hybrid-kinetic simulations of ion heating in Alfv\'enic turbulence. The Astrophysical Journal, 2019, 879(1): 53.

\bibitem{cancellation} Stanier A, Chacon L, Le A. A cancellation problem in hybrid particle-in-cell schemes due to finite particle size. Journal of Computational Physics, 2020, 420: 109705.

\bibitem{ha1} Crouseilles N, Einkemmer L, Faou E. Hamiltonian splitting for the Vlasov--Maxwell equations. Journal of Computational Physics, 2015, 283: 224-240.

\bibitem {ha2} Qin H, He Y, Zhang R, Liu J, Xiao J, Wang Y. Comment on "Hamiltonian splitting for the Vlasov--Maxwell equations". Journal of Computational Physics, 2015, 297: 721-723.

\bibitem{ha3} He Y, Qin H, Sun Y, Xiao J, Zhang R, Liu J. Hamiltonian time integrators for Vlasov--Maxwell equations. Physics of Plasmas, 2015, 22(12).

\bibitem{Leimkuhler} Bond S D, Leimkuhler B J, Laird B B. The Nos\'e--Poincar\'e method for constant temperature molecular dynamics. Journal of Computational Physics, 1999, 151(1): 114-134.

\bibitem{xuliwei} Xu L, Guyenne P. Numerical simulation of three-dimensional nonlinear water waves. Journal of Computational Physics, 2009, 228(22): 8446-8466.







\end{thebibliography}
\end{document}